\RequirePackage{fix-cm}
\documentclass[smallcondensed]{svjour3}     
\smartqed  
\usepackage{ulem,amsmath, amsfonts,amssymb,graphicx,color,rotating,url,tikz, pgffor, pgfgantt, xcolor,xkeyval,hyperref, listings, mathtools, etoolbox,bm, aliascnt,epstopdf, float,algorithm, booktabs, enumitem, multirow}
 \usetikzlibrary{arrows, automata,shapes, decorations}
\usepackage[noend]{algpseudocode}
\usepackage[numbers]{natbib}

\usepackage[nolist]{acronym}
\begin{acronym}
\acro{MP}{Myopic Policy}
\acro{RSF}{Rollout of Stationary Future}
\acro{RRO}{Rollout for Relocation Only}
\acro{LSF}{Lookahead of Stationary Future}
\acro{DNF}{Dynamic No-Flexibility}
\acro{MNF}{Myopic No-Flexibility}
\end{acronym}

\newlist{steps}{enumerate}{1}
\setlist[steps, 1]{label = \textbf{Step \arabic*:}}

\begin{document}

\title{Managing mobile production-inventory systems influenced by a modulation process
}

\titlerunning{Managing mobile production-inventory systems influenced by a modulation process}        

\author{Satya S. Malladi       \and
        Alan L. Erera          \and
        Chelsea C. White III 
}


\institute{S.S. Malladi \at
              Kantar Analytics Practice, India\\
              \email{sarvanisatya@gmail.com}           
           \and 
           A.L. Erera \at
           School of Industrial and Systems Engineering, Georgia Institute of Technology,  USA \\
           \and
           C.C. White III  \at
           School of Industrial and Systems Engineering, Georgia Institute of Technology,  USA
}


\maketitle

\begin{abstract}
The objective of this paper is to investigate the potential added value of being able to relocate production capacity, relative to fixed production capacity, in a network of multiple, geographically distributed manufacturing sites. There is a growing number of examples of production capacity that can be geographically relocated with a modest amount of effort; e.g., 3D printers, bioreactors for cell and gene manufacturing, and modular units for pharmaceutical intermediates. Such a capability shows promise for enabling the fast fulfillment of a distributed network with a reduction in the total inventory and total production capacity of a distributed network with fixed production capacity without sacrificing customer service levels or total system resilience.  Allowing also for transshipment, we model a production-inventory system with L production sites and Y units of relocatable production capacity, develop efficient and effective heuristic solution methods for dynamic relocation and multi-location inventory control, and analyze the potential added value and implementation challenges of being able to relocate production capacity. We describe the (L, Y) problem as a problem of sequential decision making under uncertainty to determine transshipment, mobile production capacity relocation, and replenishment decisions at each decision epoch. To enhance model realism, we use a partially observed stochastic process, the modulation process, to model the exogenous and partially observable forces (e.g., the macro-economy) that affect demand.  We then model the (L, Y) problem as a partially observed Markov decision process. Due to the considerable computational challenges of solving this model exactly, we propose two efficient, high quality heuristics. We show for an instance set with five locations that production capacity mobility and transshipment, relative to the fixed production capacity case, can improve systems performance by as much as 41\% on average over the no-flexibility case and that production capacity mobility can yield as much as 10\% more savings compared to when only transshipment is permitted. 

\end{abstract}
\DeclarePairedDelimiterX{\abs}[1]{\lvert}{\rvert}{\ifblank{#1}{{}\cdot{}}{#1}}
\DeclarePairedDelimiter{\ceil}{\lceil}{\rceil}
\DeclarePairedDelimiter{\floor}{\lfloor}{\rfloor}
\section{Introduction}
The aim of this paper is to better understand the impact of being able to relocate easily relocatable production capacity in a distributed production-inventory system subject to time varying demand uncertainty. In achieving this aim, we develop efficient methods that can be useful in addressing such questions as: (i) When, how much, and to where inventory should be transshipped and/or transportable production capacity should be relocated? (ii) How should replenishment decisions be made in coordination with this capability to transship inventory and/or relocate production capacity? We believe that the ability to relocate production capacity can enable the fast fulfillment of a distributed network with a reduction in the total inventory and total production capacity of a distributed network with fixed production capacity without sacrificing customer service levels or total system resilience, and the results of this paper investigate this belief.  Although we recognize that seamlessly operating mobile production capacity at various locations depends on agile logistics, for modeling and computational simplicity, we will not consider the logistics implications of production capacity relocation in this paper. 

In examining the potential value added of relocatable production capacity, we investigate a multi-period, distributed production-inventory system under stochastic demand that allows backlogging, assumes instantaneous replenishment, and has the capability to relocate transportable production units and/or transship inventory between locations. Historically, transshipment has been a tool to reposition inventory in order to improve supply chain performance. We now add the capability of repositioning production capacity to further aid in improving the performance of a supply chain. Transportable production units, which we refer to as modules, have recently generated significant interest in manufacturing \cite{GeekWire2018,BayerTechnologyServicesGMBH2014,Pfizer2015,MIT2016}; we remark that manufacturing and/or storing the final, or near-final product close to demand can enable fast fulfillment. Relatedly, we further remark that relocatable storage capacity has also become of interest recently to the parcel express industry, as described in  \cite{Verlinde2014} and elsewhere.

We model this problem as a specially structured, large-scale, partially observed Markov decision process (POMDP) in order to determine replenishment decisions, when to transship and/or relocate production capacity, and hence determine the value of having the capability to transship inventory and relocate production capacity. Our model of data-driven demand forecasting assumes the existence of a stochastic process, the modulation process, that affects demand. The modulation process is governed by an action independent Markov chain and is partially observed by the demand process and another process, the additional observation data (AOD) process. The modulation process can model exogenous factors, such as current macro-economic conditions, the weather, and seasonal effects that can affect the demand process. The AOD process models all available data beyond demand data, e.g., interest rates, unemployment rates, consumer price indices, that might be useful in improving demand forecast accuracy. Including the AOD process in our model reflects the facts that supply chains are becoming increasingly data driven to support improved supply chain performance and industry has major interest in using increased data availability velocity, volume, and variety to improve demand forecasts.  We will show that the current belief function of the modulation process can significantly influence the current demand forecast. The problem objective is to minimize the expected total discounted cost criterion composed of backorder, holding, transshipment, and module relocation costs.

A complexity analysis indicates the need to develop good, tractable heuristics (i.e., sub-optimal policies) for solution determination. We approach the development of heuristics in two ways: a centralized approach and a decentralized approach. We investigate the quality and the computational characteristics of the heuristics developed. A desirable, but not guaranteed, feature of a (sub-optimal) heuristic is that the heuristic will improve the system’s performance with improved observation accuracy. We present a preliminary numerical study that indicates the heuristics under consideration share this feature with high likelihood. With regard to the value of being able to relocate production capacity, we analyze the value of transshipment without module relocation, module relocation without transshipment, and both transshipment and module relocation.  We find that in certain cases, the value added due to resource mobility can be significant, indicating the potential importance of resource mobility in the design of next generation supply chains.
  
More specifically, we consider a distributed production - inventory system with $L$ locations and $Y$ transportable production modules. None, one, or more up to a maximum number of modules can be located at each of the locations. At each decision epoch, we assume the (centralized) decision maker (DM) knows the current demand forecast, inventory level, and production capacity at each location. This production capacity is made up of fixed capacity and transportable capacity. The DM decides how the current inventory and transportable production capacity should be relocated. We assume these relocations occur instantaneously. Once the inventory and transportable production capacity have been relocated, the DM  determines the replenishment decisions at each location based on current demand forecasts, the new inventory levels, and the new production capacities at the locations. Replenishment is instantaneous. Once replenishment is complete, demands at the locations are realized. Based on these realizations and possibly other data, the demand forecast is updated just before the next decision  epoch. 
  
We remark that how often modules would be relocated and the lead time for the relocation would in reality be situation dependent. For 3D printers, bioreactors, and smart locker modules (a form of mobile storage capacity), relocation decisions might be made monthly and take a day or two, with other decisions (e.g., reagent replenishment for the bioreactors) made hourly, daily, or weekly, where reagent replenishment may have a two-week lead time. When disruptions occur, these decisions may be event driven. Relocation decisions for modular production units for pharmaceutical intermediates might be made quarterly and require several weeks of lead time if the unit is to be moved across a national border. However, if the relocatable module is a smart locker or 3D printer on a truck or
trailer  \cite{GeekWire2018, Verlinde2014}, then relocation decisions might be made once, possibly twice,
daily and require less than an hour or two of lead time. Since the aim of the research presented in this paper is to understand the potential impact of being able to relocate production capacity, for modeling and computational simplicity, we will
assume both replenishment and relocation decisions are made instantaneously at each decision epoch with full knowledge that a specific application would require a more realistic, and a more computationally challenging, model.  

\subsection{\label{lit}Literature Review}
The problem considered in this paper involves inventory transshipment, mobile production capacity relocation, fixed production capacity of each single location production facility, and a centralized controller determining transshipment, module relocation, and replenishment decisions. Numerous innovative developments in manufacturing, such as containerized production for pharmaceutical manufacturing processes \cite{BayerTechnologyServicesGMBH2014,Pfizer2015,MIT2016} and on-demand mobile production \cite{GeekWire2018} necessitate the planning of logistics for flexible production and inventory systems that are characterized by resource mobility, interconnectivity, sharing, and decentralization \cite{Marcotte2016IntroducingProduction}. 
Malladi\textit{ et al.}~\cite{Malladi2018a} investigate the dynamic mobile production and inventory problem without the option of inventory transshipment under stationary and independent demands and have proposed heuristic approaches to solve the problem. A value addition of more than 10\% over in-the-ground production systems was determined for a system of  twenty locations.  W\"{o}rsd\"{o}rfer \textit{et al.}~\cite{Worsdorfer2017} present a real options pricing based method of evaluating the value added by mobile containerized production systems. Other research that address the operational logistics of mobile facilities can be found in \cite{Halper2011,Qiu2013}. The problem of managing mobile production capacity under deterministic demands may be viewed as a dynamic facility location problem with multiple facilities at the each location that may be opened and closed \cite{ghiani2002, melo2006, jena2015, worsdorfer2017b}. However, inventory is generally not managed jointly with capacity allocation in these problems. Solving an expanded mixed integer program, which is often the solution approach proposed in literature, will not be tractable under uncertainty and inventory control in tandem. Additionally, a mixed integer programming approach may not even be able to incorporate complex demand processes with a large number of potential demand outcomes, such as the one addressed in the current paper.

Regarding multi-location inventory management with transshipment, Karmarkar \cite{Karmarkar1979,Karmarkar1981} considers the multi-location inventory control problem over a single period and multiple periods, respectively, under uncertain demands. It is proved that when the inventory addition and subtraction matrix has a Leontief structure, there exists a base stock policy that is optimal when attainable. In \cite{Karmarkar1987}, a restricted Lagrangean dual -based lower bound and a dual relaxation based upper bound on the optimal cost of the multi-location problem are presented. The upper bound assumes the post ordering and shipment inventory position does not fall below the initial inventory position.  Rudi \textit{et al.}~\cite{Rudi2001AMaking} indicate that localized transshipment strategies are outperformed by centralized strategies. Axsater \textit{ et al.}\cite{Axsater2002} propose heuristics for a problem that considers inventory held at a warehouse and allocated for distribution to various locations in a centralized fashion. Herer \textit{ et al.} \cite{Herer2006} prove the optimality of order-up-to policies at each location in a multi-location inventory control system with reactive transshipment for a long-run average cost criterion and present a heuristic for computation. The authors consider only replenishment decisions that result in non-negative inventory positions post replenishment at each location. Lien \textit{ et al.} \cite{Lien2011} present a comparison of chain and group configurations of transshipment network design building on the ideas of manufacturing process flexibility \cite{Jordan1995} and restricted connectivity in a transshipment network \cite{Herer2002}.
Wee \textit{ et al.} \cite{Wee2005} consider a multi-retailer, one warehouse framework that allows reactive transshipment either from the warehouse to the retailers and/or between retailers. The authors prove that it is optimal to adopt either retailer only, warehouse only,  retailer first, or warehouse first protocols, when considering transshipment. Various cost parameter thresholds based intervals are presented to indicate the system that is optimal in each regime. 

The literature on the single location inventory problem is vast and varied \cite{Malladi2018, Katehakis2015, cheung15, godfrey01, Bernstein2006}. We consider the data-driven online learning demand model presented in \cite{Malladi2018} and adopt it for the multi-location problem in this paper. Malladi \textit{ et al.} \cite{Malladi2018} analyze a single location, infinite capacity inventory control problem with demand and additional observation data influenced solely by a Markov \textit{modulation} process. The modulation process is intended to model forces that may be partially observed, influence the demand process, but are not affected by actions taken by the DM (e.g., the macro-economy, air currents, tides).  Demand realizations and other data (e.g., housing starts, consumer spending) represent observations of the modulation process. What is known to the DM about the modulation process is provided by the belief function, which is updated with new data using Bayes’ Rule. A base stock policy, having a base stock level dependent on the belief function, is proved to be optimal for the infinite horizon problem when an attainability assumption holds.  The modulation process can be used to model the correlation between demands at different locations.

We consider approximate dynamic programming approaches that do not rely on maintaining the cost function's lookup table over the entire horizon to find good heuristic solutions to the multi-location mobile capacity and inventory control problem \cite{Powell2012,Ryzhov2012,Secomandi2001,Goodson2017,adpbook, Apostolos1997}. In particular, we are interested in rollout based heuristics which are known to perform well on dynamic systems with stochasticity as suggested in \cite{Secomandi2001} for solving the vehicle routing problem with stochastic demands. Goodson \textit{ et al.} \cite{Goodson2017} provide a systematic classification-aimed analysis of rollout policies. Additionally, the literature suggests that centralized control is expected to perform better than decentralized control from a solution-quality perspective; however, there is an inherent tradeoff between solution quality and computational expense \cite{Kouvelis1997TheSystem,Bernstein2005,Bernstein2006}. In the current paper, we propose and analyze a decentralized control policy that performs comparably with a more computationally intensive centralized control policy. Needed foundational results can be found in the appendices.  
\subsection{Paper Outline}
This paper is organized as follows. In \autoref{prel_res}, we state the problem, model it as a POMDP, present several preliminary results for the POMDP, and examine the tractability challenges of this model.  These challenges indicate the need for heuristic approaches. In \autoref{bounds}, we develop an approximation of the value function for the general $L$ production facility model, based on the value function of the least computationally demanding, single production facility problem. 
We also discuss the challenges of solving the $L=1$ case.  We then determine two approximations for the $L=1$ problem in \autoref{approx_single}.  In \autoref{heurist}, we present five heuristics for solving the general problem, based on these two approximations.  \autoref{comput} then presents the results of a computational study of these five heuristics. We observe performance improvement when production capacity is mobile as high as 26\% in some instances, relative to systems with no mobility, irrespective of the presence of transshipment flexibility. Also, we note that non-stationary modeling of demand when demand is non-stationary, rather than using a stationary approximation, can result in as much as a 6\% increase in performance and that complete observability of the modulation process can increase the value addition of mobility by 5\% to 27\% on the instances considered. Additionally, we infer that although joint control results in slightly lower costs, decentralized control heuristics require significantly less computational time.  Conclusions are presented in \autoref{concl}. Appendix \ref{struc_L1} contains significant foundational results for the $L=1$ case that are needed for the analysis of the heuristics, Appendix \ref{proof_prop7} presents the proof of a key result, Appendix \ref{ex_hr} presents a computationally useful heuristic, and Appendix \ref{ex_tables} presents additional numerical results that complement results in \autoref{comput_res}.

\section{Problem Statement and Preliminary Results \label{prel_res}}
We now define the general $L$ location, $Y$ module problem statement in \autoref{pr1} and present the POMDP model of this problem and general results for the model in \autoref{pr2}.  In \autoref{bounds}, we will examine the simplest case ($L = 1$, $Y = 0$) and use its solution as the basis for the development of heuristics for the general problem.  

\subsection{Problem Statement\label{pr1}}
Consider a distributed production-inventory system with $L$ locations and $Y$ portable manufacturing modules. At each decision epoch $t$ we assume the decision-maker (DM) knows:
\begin{itemize}
\item $\bm{s}(t) = \{s_l(t), l=1, \dots, L\}$, where $s_l(t)$ is the inventory level at location $l$, 
\item $\bm{u}(t) = \{u_l(t), l=1,\dots, L\}$, where $u_l(t) \in \{0,1,\dots, Y_l'\}$ is the number of modules positioned at location $l$ and $Y_l'$ is the maximum number of modules that location el can support (e.g., due to space limitations), 
\item $\mathcal{I}(t) = \{\bm{d}(t), \dots, \bm{d}(1), \bm{z}(t), \dots, \bm{z}(1), \bm{x}(0)\}$, where:
\begin{itemize}\setlength{\itemindent}{-.1in}
\item $d_l(t)$ is the demand realized during period $(t-1, t)$ that location $l$ is required to fulfill (or back order) and $\bm{d}(t) =\{d_l(t), l=1,\dots, L\}$
\item $\bm{z}(t)$ represents data, in addition to the realization of demand, that might be of use to the DM,
\item $\bm{x}(0)$ is an a priori  probability vector defined below.
\end{itemize} 
\end{itemize}
We assume the \textit{demand} process $\{\bm{d}(t), t=1, 2, \dots \}$ and the \textit{additional observation data (AOD)} process $\{\bm{z}(t), t=1, 2, \dots\}$ are linked to the \textit{modulation} process $\{{\mu}(t), t=0,1, \dots \}$ through the given conditional probability $P(\bm{d}(t+1), \bm{z}(t+1), {\mu}(t+1) \mid {\mu}(t) )$, where $\bm{x}(0) = \{x_i(0),| i \in \{1,\dots,N\}\}$  where $x_i(0) = P({\mu}(0) = \mu_i)$ for each of $N$ modulation states.  A discussion of this general description of data-driven demand and learning and how it generalizes and extends the Markov-modulated demand and Bayesian updating literatures can be found in \cite{Malladi2018}. 

The chronology of events within period $(t, t+1)$ is as follows:
\begin{steps}[wide, labelwidth=!, labelindent=0pt]
\item Given $\mathcal{I}(t)$, $\bm{s}(t)$, and $\bm{u}(t)$, the DM relocates inventory and modules to reach the post-movement state $(\bm{s}'(t),\bm{u}'(t))$, where we assume $\sum_{l=1}^L s_l'(t) = \sum_{l=1}^L s_l(t)$ and $\sum_{l=1}^L u'_l(t) = \sum_{l=1}^L u_l(t)$. 
Necessarily, $- (s_l(t))^+ \leq $ $ \Delta^S_l(t) \leq $ $ \sum_{k=1,k\neq l}^L$ $  (s_k(t))^+ $ for each location $l$, where $\Delta^S_l(t)$ is the amount of inventory relocated to location $l$. Thus, $s_l'(t) = s_l(t)+\Delta^S_l(t)$ for all $l$ and hence $\bm{s}'(t) = \bm{s}(t)+\bm{\Delta^S}(t)$, where $\bm{\Delta^S}(t) = \{ \Delta^S_l(t), l = 1, \dots, L\}$. The decision variables are $\bm{\Delta^S}(t)$ and $\bm{u}'(t)$ for Step 1. 
\item Given $\mathcal{I}(t)$, $\bm{s}'(t)$, and $\bm{u}'(t)$, the DM determines $\bm{q}(t) = \{q_l(t), l=1,\dots, L\}$, where $q_l(t)$ is the replenishment decision at location $l$. 
Necessarily, $0\leq q_l(t)$ $ \leq$  $U_l+u_l'(t) G$, where $U_l$ is the fixed amount of capacity at location $l$ and $G$ is the capacity of each module. Let $y_l(t) = s_l'(t) + q_l(t)$, the inventory level at location $l$ after inventory and module relocation and replenishment but before demand realization, and assume $\bm{y}(t) = \{y_l(t), l=1, \dots, L\}$. The decision variables are therefore $\bm{q}(t)$, or equivalently $\bm{y}(t)$, for Step 2, where necessarily, $s_l'(t) \leq y_l(t) \leq s'_l(t)+U_l + u_l'(t) G$ for all $l$. 
\item The realizations of the random variables $\bm{d}(t+1)$ and $\bm{z}(t+1)$ become known and unfulfilled demands are backordered, $\mathcal{I}(t+1) = \{ \bm{d}(t+1), \bm{z}(t+1), \mathcal{I}(t)\}$, $\bm{s}(t+1) = \bm{y}(t) - \bm{d}(t+1)$, and $\bm{u}(t+1)= \bm{u}'(t)$.  
\item $t=t+1$. 
\end{steps}
We assume that for location $l$, $c_l(y_l(t),d_l(t+1)) = b_l (d_l(t+1)-y_l(t))^+ +h_l(y_l(t)-d_l(t+1))^+ \geq 0$ is the single period cost accrued between $t$ and $t+1$, where $b_l$ and $h_l$ are respectively the backorder and holding cost per unit per period and for all $d_l$, $c_l(y_l,d_l)$ is convex in $y_l$ and $\lim_{\abs{y} \to \infty} c_l(y, d_l) = \infty$. 

We assume that the modulation and the observation state spaces are finite and that for each location, the demand state space is finite and the inventory state space $\{\dots, -1,0,1, \dots\}$ is countable.

Let the single period $(t, t+1)$ cost be:
\begin{multline*}
\vspace{-2ex} 
\sum_{l=1}^L \left(K^{S+}_l (\Delta^S_l(t))^+ + K^{S-}_l (-\Delta^S_l(t))^+ \right) 
\\+ K^M\sum_{l=1}^L \abs{u_l(t)-u_l'(t)}/2 + \sum_{l=1}^L c_l \big(y_l(t), d_l(t+1) \big),
\end{multline*}
where $K^{S+}_l$ ($K^{S-}_l$) is the cost of moving a unit of inventory to (from) location $l$, and $K^M$ is the cost of moving a module from one location to another. A \textit{feasible} policy determines $\big(\bm{q}(t), \bm{\Delta^S}(t), \bm{u}'(t) \big)$ based on $\mathcal{I}(t), \bm{s}(t),$ and $\bm{u}(t)$ for all $t$. 

The problem criterion is the expected total discounted cost over the infinite horizon, where $\beta \in [0,1)$ is the discount factor. The problem is to determine a feasible policy that minimizes the criterion with respect to the set of all feasible policies. 

We remark that we have used the model of transshipment cost described above due to its modeling simplicity and note that for some transshipment problems (e.g., between retail stores within an urban area) this model might be reasonably suitable.  However, in general the cost of transshipment will be different for different origin-destination pairs, and hence a specific real-world example of a distributed production-inventory problem may require a more realistic model of transshipment costs.  
\subsection{POMDP Model and General Results\label{pr2}}
This problem can be recast as a partially observed POMDP as follows. Results in \cite{Smallwood1973} and \cite{Sondik1978TheCosts} imply that $(\bm{x}(t),\bm{s}(t), \bm{u}(t))$ is a sufficient statistic, where the \textit{belief} function $\bm{x}(t) = \{x_i(t), \forall i=1,\dots,N\}$ is such that $x_i(t) = P(\mu(t) = \mu_i \mid \mathcal{I}(t))$ and $\bm{x}(t) \in X = \{\bm{x}\geq 0: \sum_{i=1}^N x_i =1 \}$. Let 
\begin{align*}
&P_{ij}(\bm{d},\bm{z}) 
\\&= P(\bm{d}(t+1) = \bm{d}, \bm{z}(t+1) = \bm{z}, \mu(t+1)=j \mid \mu(t) = i) 
\\ & \forall i,j \ \in \{1,\dots, N\},
\\ &\sigma(\bm{d},\bm{z},\bm{x}) = \bm{x}P(\bm{d},\bm{z})\underline{1} = \sum_{i=1}^N  x_i \sum_{j=1}^N P_{ij}(\bm{d},\bm{z}),
 \\ & \bm{\lambda}(\bm{d},\bm{z},\bm{x}) = \{ \lambda_j(\bm{d},\bm{z},\bm{x}), \forall j = 1,\dots, N \} 
 \\ & = \bm{x}P(\bm{d},\bm{z}) / \sigma(\bm{d},\bm{z},\bm{x}), \ \sigma(\bm{d},\bm{z},\bm{x}) \neq 0, \text{ and }
\\ & {\mathcal{L}}(\bm{x},\bm{y}) = E[{c}(\bm{y},\bm{d})] 
\\ & = \sum_{\bm{d},\bm{z}} \sigma(\bm{d},\bm{z},\bm{x}) {c}(\bm{y},\bm{d}), \ \ \  {c}(\bm{y},\bm{d}) = \sum_{l=1}^L c_l(y_l,d_l). 
\end{align*}
Thus, if $\bm{x}$ is the prior belief function, then $\bm{\lambda}(\bm{d},\bm{z},\bm{x})$ is the posterior belief function, given realizations $(\bm{d},\bm{z})$, and $\sigma(\bm{d},\bm{z},\bm{x})$ is the probability that $(\bm{d},\bm{z})$ will be the demand and observation realizations, given prior $\bm{x}$. 
Define the operator ${H}$ as follows:
\begin{eqnarray} \label{eq1} [{H}{v}](\bm{x},\bm{s},\bm{u}) = \min_{\bm{\Delta^S}, \bm{u}', \bm{y}} \{ {\mathcal{G}}(\bm{x}, \bm{u}, \bm{y},\bm{\Delta^S}, \bm{u'}, {v}) \},  \end{eqnarray}
\begin{multline*}
\text{where }
{\mathcal{G}}(\bm{x}, \bm{u}, \bm{y},\bm{\Delta^S}, \bm{u'}, {v})
\\= \sum_{l=1}^L \left( K^{S+}_l (\Delta^S_l)^+ 
 + K^{S-}_l (-\Delta^S_l)^+ \right)
 \\ + K^M\sum_{l=1}^L \abs{u_l -u_l'}/2 + {\mathcal{L}}(\bm{x},\bm{y}) 
\\ + \beta \sum_{\bm{d},\bm{z}} \sigma(\bm{d},\bm{z},\bm{x}) {v}\big(\bm{\lambda}(\bm{d},\bm{z},\bm{x}),\bm{y}-\bm{d}, \bm{u}' \big),
\end{multline*} 
and where the minimization is with respect to 
\begin{eqnarray*} && \sum_{l=1}^L u_l' = Y,
 \\  && 0\leq u_l' \leq Y_l', \ \forall \ l \in \{1, \dots, L \},
 \\ && \sum_l \Delta^S_l =0,
  \\ && - (s_l)^+ \leq \Delta^S_l \leq \sum_{k\neq l} (s_k)^+ , \ \forall \ l \in \{1, \dots, L \},
\\ && (s_l + \Delta^S_l)\leq y_l \leq (s_l + \Delta^S_l)+U_l +u_l'G, 
\\ && \forall \ l \in \{1, \dots, L \}, \text{ and }
 \\ && u_l',\ \Delta^S_l, \ y_l \in \mathbb{Z}, \ \forall \ l \in \{1, \dots, L \}.
 \end{eqnarray*} 
Results in Puterman \textit{ et al.} \cite{Puterman1994} guarantee that there exists a unique ${v}^*$ such that ${v}^* = {H}{v}^*$ and that this fixed point is the minimum expected total discounted cost over the infinite horizon. Further, a  policy that causes the minimum in  \eqref{eq1} to be attained is an optimal policy and is decision epoch invariant. For any given bounded function ${v}_0$, let $\{{v}_n\}$ be such that ${v}_{n+1} = {H}{v}_n$. Then, $\lim_{n\to \infty} \abs{\abs{{v}^*-{v}_n}} = 0$, where $\abs{\abs{.}}$ is the sup-norm. 

Results in \cite{Smallwood1973} guarantee that ${v}_n(\bm{x},\bm{s},\bm{u})$ is piecewise linear and concave in $\bm{x}$ for fixed $(\bm{s},\bm{u})$ for all $n$, assuming ${v}_0(\bm{x},\bm{s},\bm{u})$ is also piecewise linear and concave in $\bm{x}$ for fixed $(\bm{s},\bm{u})$. In the limit, ${v}^*(\bm{x},\bm{s},\bm{u})$ may no longer be piecewise linear in $\bm{x}$ for fixed $(\bm{s},\bm{u})$;  however, concavity will be preserved. The cardinality of the state space of this POMDP is infinite since the belief vector belongs to a continuous $N$-dimensional real space. 
Even when dealing with the i.i.d. case, since the number of inventory and module count combinations is exponential in the number of locations $L$, solving this POMDP exactly becomes intractable for even relatively small values of $L$.
Therefore, we seek good sub-optimal approaches that significantly reduce this computational burden. 
\section{Bounds and Approximate Value Function Based on $L=1$ Case \label{bounds}}
Throughout this paper, we will  base the development of heuristics on the most tractable problem, the single location inventory control problem, i.e., the $L=1, Y=0$ case. Solving each of the $L$ local replenishment problems for the i.i.d. case requires $\abs{S_l}^2  \abs{A_l}$ multiplications per successive approximation iteration, and $L$ of these are required. For $L=10$ and $\abs{S_l}$ = $ \abs{A_l} = 50$, $L\abs{S_l}^2\abs{A_l}$ is on the order of $10^{5}$, which is a large but computationally manageable problem.  The operator ${H}$ simplifies to $H^F_l$ for location $l$ with \textit{fixed} capacity, where 
\begin{eqnarray} \label{eq_1a}
 && [H^F_l v^F_l]  (\bm{x},s_l, u_l) =\min \big \{\mathcal{G}^F_l(\bm{x},u_l , y_l,v^F_l) \big\}, 
\\ && \nonumber \mathcal{G}^F_l(\bm{x},u_l, y_l, v^F_l) =\mathcal{L}^F_l (\bm{x}, y_l) 
\\ &&  \nonumber + \beta \sum_{\bm{d},\bm{z}} \sigma(\bm{d},\bm{z},\bm{x}) v^F_l(\bm{\lambda}(\bm{d},\bm{z},\bm{x}), y_l-d_l, u_l), 
\\ &&  \nonumber \text{where } \mathcal{L}^F_l (\bm{x}, y_l) = \sum_{\bm{d},\bm{z}} \sigma(\bm{d}, \bm{z}, \bm{x}) c_l(y_l, d_l), 
\\ &&  \nonumber \forall \ l \in \{1, \dots, L \}, 
\end{eqnarray}
and where the minimization in \eqref{eq_1a} is with respect to $s_l \leq y_l \leq s_l+U_l + u_l G$. \autoref{prop3} in Appendix  \ref{struc_L1} guarantees that the fixed point of $H^F_l$, $v^F_l$, is non-decreasing in capacity for fixed $(\bm{x},s_l)$. This monotonicity result implies
$$ \sum_{l=1}^L v^F_l(\bm{x},s_l, Y_l') \leq {v}(\bm{x},\bm{s},\bm{u}).$$
At this point, it is important to note that the arguments for the function $\mathcal{G}$ in \eqref{eq1} are different from the arguments for the function $\mathcal{G}^F_l$ in \eqref{eq_1a}. The arguments of $\mathcal{G}$ contain the additional terms $\bm{\Delta^S}$ and $\bm{u'}$, which are the relocation decision variables in \eqref{eq1}. Implicit in these terms being absent in the arguments of $\mathcal{G}^F_l$ is the assumption that for the single location case, transshipment and/or module relocation are assumed not to occur in the future. 
Hence,
$${v}(\bm{x},\bm{s},\bm{u}) \leq \sum_{l=1}^L v^F_l(\bm{x},s_l, u_l).$$
Thus, the solutions of the local replenishment problems provide  upper and lower bounds on the cost function of the initial problem. 

We now approximate the optimal cost-to-go function of the POMDP presented in \eqref{eq1}. 
Let $\theta \in [0,1]$, be such that 
 \begin{eqnarray*} {v}_l^{F,\theta} (\bm{x}, s_l, u_l) &=& (1-\theta) v^F_l (\bm{x},s_l, Y_l')+ \theta v^F_l (\bm{x},s_l, u_l) 
 \\&& \forall \ l \in \{1, \dots, L \}  \text{ and  } 
 \\ \tilde{v}^{\theta}(\bm{x},\bm{s},\bm{u}) &=& \sum_{l=1}^L v_l^{F,\theta} (\bm{x}, s_l, u_l). 
 \end{eqnarray*}
Hence, $\tilde{v}^{\theta}(\bm{x},\bm{s},\bm{u})$ is an approximation of ${v}(\bm{x},\bm{s},\bm{u})$ that relies solely on the solution of the single location ($L=1, Y=0$) problem. Then, 
 \begin{align}
\nonumber [{H}\tilde{v}^{\theta}] &(\bm{x},\bm{s},\bm{u}) = \min_{\bm{\Delta^S}, \bm{u}'} \bigg\{ \sum_{l=1}^L ( K^{S+}_l (\Delta^S_l)^+ 
\\ &  + K^{S-}_l (-\Delta^S_l)^+ ) 
+ K^M\sum_{l=1}^L \abs{u_l- u_l'}/2  \nonumber
 \\ &
 + \min_{\bm{y}} \{ \mathcal{L}(\bm{x},\bm{y}) \nonumber
\\& +\beta \sum_{\bm{d},\bm{z}} \sigma(\bm{d},\bm{z},\bm{x}) \tilde{v}^{\theta}(\bm{\lambda}(\bm{d},\bm{z},\bm{x}), \bm{y}-\bm{d}, \bm{u}')\}\bigg\}.\label{blend_hr}
 \end{align}
 
 In \eqref{blend_hr}, the inner minimization is over all $y_l$ such that $s_l+\Delta^S_l \leq y_l \leq s_l +\Delta^S_l +U_l  + u_l' G$. We note that the coefficient $\theta$ may be location-dependent depending on the nature of the instances.  

\section{Approximating the Value Function for the Single Location Problem \label{approx_single}}
Two of the heuristics presented in \autoref{heurist} make use of the value function of the single location problem.  We remark that when $L=1$, the (local) decision maker assumes there will be no inventory and/or module relocation in the future and does not attempt to coordinate its decisions with either the controller determining the inventory and/or module relocation decisions or the replenishment decision makers at the other locations. Foundational results for the $L=1$ problem presented in Appendix \ref{struc_L1} imply that there exists an optimal replenishment policy that is a base stock policy, the optimal base stock value is non-increasing in capacity, but an optimal base stock policy is myopic only when production capacity is sufficiently large (\autoref{prop5} of \autoref{struc_L1}).  Computational complexity and the likelihood of intractability for the case where demand is not i.i.d., even for the $L = 1$ problem, increases substantially when an optimal policy is not myopic.  For computational reasons, we now present two approximations of the value function for the single location problem, the static belief function approximation and the piecewise linear approximation based on the convexity of the value function in $s$ and $u$.  
\subsection{The Static Belief Function Approximation \label{stat_L1}}
Assume that $\bm{x}(t+1) = \bm{x}(t)$ for all $t$.  Then the $L = 1$ operator becomes 
\begin{align}
 \big[ \hat{H}^F_l \hat{v}^F_l \big]& (\bm{x},s_l,u_l)  \nonumber
 \\ & = \min_{s_l\leq y_l\leq s_l+U_l+G u_l} \bigg\{  \sum_{d_l} \sum_{i=1}^N x_i \text{Pr}(d_l\mid i) \bigg [ c_l(y_l,d_l) \nonumber
\\ & + \beta \hat{v}^F_l (\bm{x},y_l-d_l, u_l) \bigg ] \bigg\} \label{stat_ap}
\end{align}
which for given $\bm{x}$ and $u_l$, requires essentially the same number of operations per successive approximations step as required in the i.i.d. case. Since the case where $\bm{x}(t+1) = \bm{x}(t)$ for all $t$ is a special case of the general problem, there exists an optimal policy that is a base stock policy, an optimal base stock level is non-increasing in capacity, and the optimal value function is non-increasing in capacity and convex in inventory level (see Appendix \ref{struc_L1}).  Thus, the resulting approximation $\hat{v}_l^F$ shares the same structural properties of $v_l^F$.
We now present a result that bounds the gap between $\hat{v}^F_l$ and $v_l^F$ that will prove useful in our computational study; proof is presented in Appendix \ref{proof_prop7}.
\begin{proposition}\label{prop7}
We have $v_l^F(\bm{x},s_l,u_l) \geq \hat{v}^F_l(\bm{x},s_l,u_l) - \rho/(1-\beta)$ for all $\bm{x}$, $s_l$, and $u_l$ for $l\in \{1, \dots, L \}$
where $\rho =\ $ $\sum_{d_l} k(d_l) c_l(\hat{y}_l,d_l)$ and $k(d_l) = \big(\max_k \text{Pr}(d_l\mid k) - \min_k \text{Pr} (d_l\mid k) \big)$.\end{proposition}
\subsection{A piecewise linear and convex approximation of the value function of $L=1$ static fixed problem \label{plc_L1}}
We use the following approximation of the optimal cost of the single location static fixed problem $\hat{v}_l^{F}$, drawing inspiration from the approximation of the cost-to-go function in the lookahead of fixed future (LAF) heuristic in \cite{Malladi2018a}: 
\begin{multline*}\hat{v}^{F}_{l} (\bm{x}(t+1), s_l(t+1), u_l(t+1)) 
\\ \approx \big(\hat{v}^{F}_{l} (\bm{x}(t+1), \overline{s}_l(t+1), u_l(t)) 
\\ + \hat{v}^{F}_{l} (\bm{x}(t+1), s_l(t), u_l(t+1)) \big)/2, \ \forall l \in \{1, \dots, L \}, 
\end{multline*}
where $\bar{s}_l(t+1) = y_l(t) - \big[E[D_l(t)]\big]$ and $\big[a\big]$ denotes the nearest integer to which $a$ is rounded. 

Since $v^F_l(\bm{x},s_l,u_l)$ is piecewise linear and convex in $s_l$ when $u_l$ is held constant and in $u_l$ when $s_l$ is held constant (from \autoref{prop1} and \autoref{prop4} in Appendix \ref{struc_L1}) and $\hat{v}^F_l(\bm{x},s_l,u_l)$ inherits these properties as it is a stationary special case, the latter  can be represented as  $\max\{ \gamma_j^l s_l + \hat{\gamma}_j^l: (\gamma_j^l, \hat{\gamma}_j^l) \in \Gamma_t^l(u_l) \}, \ \forall l \in \{1, \dots, L\}$  
 and as $\max\{ \theta_j^l u_l + \hat{\theta}_j^l: ( \theta_j^l, \hat{\theta}_j^l) \in \Theta_t^l (s_l)\}, \ \forall l \in \{1, \dots, L\}$. The set $\Gamma_t^l (u_l)$ ($\Theta_t^l (s_l)$) is the set of coefficients describing the facets of the piecewise linear and convex function $\hat{v}^F_l(\bm{x},s_l,u_l)$, when $u_l$ ($s_l$) is held constant at time $t$. 
 Thus, the following expression is the approximation:
\begin{align}
\nonumber & \hat{v}^F_{l} (\bm{x}(t+1), s_l(t+1), u_l(t+1)) 
\\ & \approx \bigg(\max\{ \gamma_j^l s_l + \hat{\gamma}_j^l: (\gamma_j^l, \hat{\gamma}_j^l) \in \Gamma_t^l(u_l(t))\} \nonumber
\\ & + \max\{ \theta_j^l u_l + \hat{\theta}_j^l: ( \theta_j^l, \hat{\theta}_j^l) \in \Theta_t^l (s_l(t)) \} \bigg)/2, \nonumber
\\ & \forall \ l \in \{1, \dots, L \}.
\end{align}
\section{Heuristics \label{heurist}}
As the cardinality of our state space is exponential in the number of locations (see \cite{Malladi2018a}), we pursue approximate dynamic programming methods \cite{bertsekas97, adpbook} to design policies instead of obtaining a representation of the entire lookup table of the optimal cost function. We first present the \acf{MP} in \autoref{myopic} to determine dynamic inventory and relocation decisions myopically, followed by two policies in \autoref{no_flex}, \acf{MNF} and \acf{DNF} ,  that do not allow inventory and module relocation. \ac{DNF} serves as our computational benchmark policy against which we compare the quality of the remaining heuristics.  We consider a class of heuristic policies known as lookahead policies, which use an approximate cost-to-go term in the  optimality equations at every decision epoch. We employ rollout policies that determine actions at every epoch by solving a forward pass of the optimality equation with the cost-to-go approximated  as the expected cost of a given policy under a specified set of conditions from the next decision epoch onward \cite{Goodson2017, Secomandi2001}. Specifically, we propose two policies which assume at every decision epoch that from the next epoch onward, mobility of production capacity and transshipment capability are not available and demand distributions remain stationary at the current belief-mixed distributions. In \autoref{jr}, we consider a rollout policy, \acf{RSF}, that determines module and inventory relocation decisions as well as production decisions at each epoch, with the described future conditions beginning from the next decision epoch. We present a policy, \acf{LSF}, that uses the same idea as \ac{RSF} but with a piecewise linear approximation of the cost-to-go term  in \autoref{laj}. In \autoref{glr}, we propose a second rollout policy, \acf{RRO}, that determines only the module and inventory relocation decisions at each epoch, with the described future conditions beginning before the production event in the current period. 
 
Here, we present additional notation that will be useful in this section. Let: $P_{ij}=\sum_{\bm{d},\bm{z}} P_{ij}(\bm{d},\bm{z}) = P(\mu(t+1) = j \mid \mu(t) = i)$ for all $i,j \in \{1, \dots, N \}$, $P = \{ P_{ij} \} $,  and $\bm{\pi}$ satisfy $\bm{\pi} = \bm{\pi} P$. Thus, $P$ is the transition matrix of the modulation process, and $\bm{\pi}$ is a stationary probability vector, which we will assume is unique in $X$ and hence has interpretation as the  distribution of the modulation process. Further, let $O_{ij}(\bm{d},\bm{z}) = P(\bm{d}(t+1) = \bm{d}, \bm{z}(t+1) = \bm{z} \mid \mu(t+1) = j, \mu(t) = i) = P_{ij}(\bm{d},\bm{z})/P_{ij}$, or equivalently, $P_{ij}(\bm{d},\bm{z}) = O_{ij}(\bm{d}, \bm{z}) P_{ij} $. Thus, $O_{ij}(\bm{d},\bm{z})$ describes the relationship between the modulation process and the demand and the AOD observations of the modulation process.
\subsection{\acf{MP}}\label{myopic}
For the \acf{MP}, the decision-maker optimizes over the one period cost function to determine relocation and replenishment decisions. At every decision epoch with current state $(\bm{x}, \bm{s}, \bm{u})$, we therefore solve the following integer program:
\begin{align}
 \nonumber \text{MP: }  \min_{\bm{\Delta^S}, \bm{u}',\bm{y}} \sum_{l=1}^L  \bigg\{ (K^{S+}_l \Delta^{S+}_l + K^{S-}_l \Delta^{S-}_l ) &
 \\  \nonumber + K^M \abs{u_l-u_l'}/2 + +   \sum_{n=1}^M \sigma(d_l^n,\bm{x}) \big[ h_l r_l^n + b_l o_l^n \big] \bigg\} &
\\ \nonumber \text{subject to } \sum_{l=1}^L u_l' = Y, &
  \\ \nonumber 0\leq u_l' \leq Y_l', \ \forall \ l \in \{1, \dots, L \}, &
   \\ \nonumber \sum_{l=1}^L \Delta^{S+}_l = \sum_{l=1}^L \Delta^{S-}_l, &
  \\ \nonumber  0\leq \Delta^{S+}_l \leq \sum_{k\neq l} (s_k)^+, \ \  0\leq \Delta^{S-}_l \leq - (s_l)^+, &
  \\  \nonumber  \forall \ \ l \in \{1, \dots, L \}, &
\\ \nonumber  (s_l + \Delta^{S+}_l -  \Delta^{S-}_l)\leq y_l \leq (s_l + \Delta^{S+}_l - \Delta^{S-}_l) &
\\ \nonumber +U_l +u_l'G, \ \forall \ l \in \{1, \dots, L \}, &
 \\ \nonumber r_l^n \geq y_l -d_l^n, \ \ o_l^n \geq d_l^n - y_l, &
 \\   \nonumber \forall \ l \in \{1, \dots, L \}, n \in \{1, \dots, M\}, &
 \\ \nonumber r_l^n, \ o_l^n \in \mathbb{Z}^+, \ u_l',\ \Delta^{S+}_l, \Delta^{S-}_l, \ y_l \in \mathbb{Z}, &
 \\  \eta_l, \ \zeta_l \in \mathbb{R}  \ \ \forall \ l \in \{1, \dots, L \}, &
 \end{align} 
  where $M$ is the number of demand outcomes at any location $l$ and where we have assumed $O_{ij}(\bm{d}, \bm{z})$ is independent of $i$ and $z$, 
  $$P(\bm{d}(t+1) \mid \mu(t+1)) = \Pi_l P(d_l(t+1) \mid \mu(t+1) ),  $$
  and $d_l^n$ is the $n$th realization of the random variable $d_l$. 
MP accounts for transshipment quantities entering and leaving each location $l$ as $\Delta^{S+}_l$ and $\Delta^{S-}_l$ respectively, the post module movement capacity count as $u_l'$, the post-replenishment inventory position as $y_l$, and the held and backlogged inventory quantities as $r_l^n$, and $o_l^n$ for the $n$th demand scenario. 
The flow balance constraints for modules and inventory are followed by the inventory accounting constraints.
We will find later that the computational quality of \ac{MP} is poor, emphasizing the need for policies that enable dynamic optimization. 
We do not use \ac{MP} as a benchmark as the computational analysis in \cite{Malladi2018a} indicates that the quality of the myopic policy is influenced by the number of locations in the system.
Thus, in the following subsection, we pursue benchmark policies that do not allow resource mobility. 
\subsection{No-Flexibility Policies}  \label{no_flex}
In this section, we present two No-Flexibility policies that provide an upper bound on the optimal solution of the $L$ location, $Y$ module problem. 
\subsubsection{\acf{MNF} Policy}
We remark that a natural and easily computed and implemented sub-optimal policy for the finite capacity $L = 1$ problem is to order either the difference between the optimal base stock value for the infinite capacity case and the current inventory level or to order the capacity of the production system, whichever is smaller.   More specifically, the local order up to level at each location $l$ is given by 
\begin{align} \label{eq_ponews}\hat{y}_l = \min \big\{\max\{s^*_l(\bm{x}),   s_l+\Delta^S_l\}, 
 s_l+\Delta^S_l+ U_l+u_l' G \big\}, \end{align} 
where $s^*(x)$ is an optimal myopic base stock level for the infinite capacity problem, as proposed by Malladi \textit{ et al.} \cite{Malladi2018}.  

The \ac{MNF} policy does not permit inventory and module relocation, assumes that local replenishment is based on the policy presented in \eqref{eq_ponews}, and assumes that the fixed, static production capacities at the locations are selected in order to minimize multi-location expected total cost with  stationary belief distribution $\pi$. We have initially considered its use as a benchmark policy owing to its performance in the single location problem and its computational simplicity; however we find that it is outperformed as an upper bound by a dynamic policy for the no-flexibility system proposed in the next subsection. 
\subsubsection{\acf{DNF} Policy: The Benchmark Policy}
The \acf{DNF} policy  is a dynamic policy of executing inventory control at each location for a fixed production module configuration, disallowing both stock and module relocations. The following integer program that accounts for the future cost must be solved at every decision epoch to implement the \ac{DNF} policy. We make use of the static belief approximation $\hat{v}_l^{F}$  of the $L=1$ subproblems' solutions (from \autoref{stat_L1}) in the future cost term. Additionally, we assume that for every $ l \in \{1, \dots, L \}$,  the local controller approximates $\bm{\hat{\lambda}}(\bm{d},\bm{z},\bm{x})$ as $\bm{\hat{\lambda}}(d_l^n,\bm{x})$  by using the data sourced locally (i.e., $d_l^n$, the demand realized at location $l$).  This assumption allows the decomposition of the cost-to-go term by location. 

We note that the \ac{DNF} policy is independent of the  coefficient $\theta$, making it a suitable benchmark.  

\begin{align}
 & \nonumber \text{DNF: }  \min \sum_{l=1}^L\sum_{q=0}^{U_l+ \tilde{u}_l G} w(l, q)  \bigg\{ \sum_{n=1}^M \sum_{i=1}^N x_i \sum_{j=1}^N P_{ij} O^l_{nj} \bigg[  h_l r_l^n 
 \\ & \nonumber + b_l o_l^n 
 +  \beta \hat{v}_l^{F}(\bm{\hat{\lambda}}(d_l^n,\bm{x}), s_l+\Delta^S+q-d^n_l, u_l+\Delta^M)  \bigg] \bigg\}, 
\\ \nonumber & \text{subject to }  
\\  \nonumber &  r_l^n \geq s_l+\sum_{\Delta^S, \Delta^M, q} w(l,\Delta^S, \Delta^M, q) \ (\Delta^S + q) -d_l^n,  
\\ \nonumber & \forall \ l\in \{1,\dots,L\}, n \in \{1,\dots,M\}, 
\\\nonumber  & o_l^n \geq d_l^n - s_l-\sum_{\Delta^S, \Delta^M, q} w(l,\Delta^S, \Delta^M, q) \ (\Delta^S + q), 
\\  & \forall \ l\in \{1,\dots,L\}, n \in \{1,\dots,M\}, \nonumber 
 \\ \nonumber &  r_l^n, \ o_l^n \in \mathbb{Z}^+ \ \forall \ n  \in \{1,\dots,M\}, 
 \\ \nonumber & w(l, q) \in \{0,1\},  \ \forall \  q \in \{0,\dots,U_l + \tilde{u}_l G \},  
 \\ &  \text{ for }
  l \in \{1,\dots,L\}, 
 \end{align}  
 where $M$ is the number of demand outcomes at any location $l$.
For this integer program, as the future cost term  is obtained from a lookup table and is a nonlinear expression, binary variables $w(l,q)$ are used to choose the production decisions at the current epoch. The constraints account for inventory flows, namely, of held ($ r_l^n$) and backordered ($o_l^n $) quantities. The value function approximation used for this policy remains relevant for the cases with flexibility as well.   
\subsection{\acf{RSF} \label{jr}}
The \acf{RSF} heuristic policy is based on the the approximation presented in  \eqref{blend_hr}. In \ac{RSF}, at each decision epoch with current state $(\bm{x}, \bm{s}, \bm{u})$, we require the integer program \ac{RSF} given below be solved. The resulting policy utilizes the same value function approximation as \ac{DNF} by assuming the local data sourcing assumption holds. 
For this integer program, since the future cost term $\hat{v}_l^{F,\theta}$ is obtained from a lookup table and is a nonlinear expression, we adopt the following formulation that uses binary variables $w(l,\Delta^S,\Delta^M,q)$ to choose the actions at the current epoch: transshipment quantity $\Delta^S$ entering location $l$, the number of modules $u$ entering location $l$, and the production quantity $a$ at location $l$. These binary variables enable suitable selection of  ${v}_l^{\theta}$ from lookup tables in the integer program:

 \begin{align}
 \nonumber & \text{\ac{RSF}:} 
 \\ \nonumber & \min \sum_{l=1}^L\sum_{\Delta^S =-(s_l)^+}^{\sum_{k\neq l} (s_k)^+} \sum_{\Delta^M = -u_l }^{Y_l'-u_l} \sum_{q=0}^{U_l +(u_l + \Delta^M) G} 
 w(l,\Delta^S, \Delta^M, q) 
 \\ & \nonumber\bigg\{ K^{S+}_l  (\Delta^S)^+ 
 + K^{S-}_l (-\Delta^S)^+ + K^M\abs{\Delta^M}/2
\\ \nonumber  &+  \sum_{n=1}^M \sum_{i=1}^N x_i \sum_{j=1}^N P_{ij} O^l_{nj} \bigg[  h_l r_l^n + b_l o_l^n 
\\ & \nonumber + \beta \hat{v}_l^{F,\theta}(\bm{\hat{\lambda}}(d_l^n,\bm{x}), s_l+\Delta^S+q-d^n_l, u_l+\Delta^M)  \bigg] \bigg\},
\\ \nonumber & \text{subject to }  
\\ \nonumber &  \sum_{l=1}^L\sum_{\Delta^S =-(s_l)^+}^{\sum_{k\neq l} (s_k)^+} \sum_{\Delta^M = -u_l }^{Y_l'-u_l} \sum_{q=0}^{U_l +(u_l + \Delta^M) G}  w(l,\Delta^S, \Delta^M, q) \ \Delta^M 
\\ \nonumber & = 0,
 \\ \nonumber &   \sum_{l=1}^L\sum_{\Delta^S =-(s_l)^+}^{\sum_{k\neq l} (s_k)^+} \sum_{\Delta^M = -u_l }^{Y_l'-u_l} \sum_{q=0}^{U_l +(u_l + \Delta^M) G}  w(l,\Delta^S, \Delta^M, q) \ \Delta^S \\\nonumber & =0,
\\ & \nonumber r_l^n \geq s_l+\sum_{\Delta^S, \Delta^M, q} w(l,\Delta^S, \Delta^M, q) \ (\Delta^S + q) -d_l^n, \\ &  \forall \ l  \in \{1,\dots,L\},\ n  \in \{1,\dots,M\}, \nonumber
\\ & o_l^n \geq d_l^n - s_l-\sum_{\Delta^S, \Delta^M, q} w(l,\Delta^S, \Delta^M, q) \ (\Delta^S + q), \nonumber
\\ &  \forall \ l  \in \{1,\dots,L\},\ n  \in \{1,\dots,M\}, \nonumber
 \\ \nonumber & r_l^n, \ o_l^n \in \mathbb{Z}^+, \ \forall l  \in \{1,\dots,L\}, \ n  \in \{1,\dots,M\}, 
 \\ \nonumber &w(l,\Delta^S, \Delta^M, q) \in \{0,1\},  
 \\ \nonumber &  \forall \ \Delta^S \in \{ -(s_l)^+, \dots, \sum_{k\neq l} (s_k)^+  \}, \text{ and }
 \\ \nonumber & \Delta^M \in \{ -u_l, \dots, Y_l'-u_l \}, 
 \\  \label{eq_jr} &  q \in \{0,\dots,U_l +(u_l + \Delta^M) G \}, \ 
  l \in \{1,\dots,L\}. 
 \end{align}  
 The first two constraints ensure the balance of module flows and transshipped inventory flows between locations. The next two sets of constraints help determine held and backordered quantities at each location.
 
 In this approach, the number of binary variables required to solve the one period problem at every epoch grows linearly in $L$ and quadratically in the total number of modules $Y$. Hence, we present a lookahead approach in Section \ref{laj} to improve the computational efficiency of the joint controller's strategy using the piecewise linear and convex approximation of $\hat{v}_l^{F,\theta}$ presented in \autoref{plc_L1} that reduces the number of binary variables used. 
\subsection{\acf{LSF} \label{laj}}
The mixed integer program \ac{LSF}, presented below, makes use of the piecewise linear and convex approximation of the single location capacitated inventory control system's cost-to-go function presented in Section \ref{plc_L1} in order to reduce the computational effort required for \ac{RSF}. Using this functional approximation of the cost-to-go function reduces the number of integer variables by $\mathcal{O}(GY^2LI)$ where G, Y, L, and I are, respectively, the capacity per module, the total number of production modules, the number of locations, and the available storage capacity at each location:
\begin{align}
 \nonumber & \text{\ac{LSF}: }
 \\ \nonumber & \min_{\bm{\Delta^S}, \bm{u}',\bm{y}} \sum_{l=1}^L  \bigg\{ (K^{S+}_l \Delta^{S+}_l + K^{S-}_l \Delta^{S-}_l ) + K^M \abs{u_l-u_l'}/2
\\ \nonumber  &+   \sum_{n=1}^M \sigma(d_l^n,\bm{x}) \bigg[ h_l r_l^n + b_l o_l^n  + \beta (\zeta_l + \eta_l)/2 \bigg] \bigg\},
\\ \nonumber & \text{subject to }  
\\ & \zeta_l \geq \gamma_j^l (y_l-\big[E[D_l(t)]\big]) + \hat{\gamma}_j^l  \ \forall \ (\gamma_j^l, \hat{\gamma}_j^l) \in \Gamma^l_{t+1}(u_l)   \nonumber
\\ & \forall \ l  \in \{1,\dots, L \}, \nonumber
\\ & \eta_l  \geq \theta_j^l u_l' + \hat{\theta}_j^l  \ \ \forall \ ( \theta_j^l, \hat{\theta}_j^l) \in \Theta^l_{t+1}(s_l)  \ \forall \ l \in \{1,\dots, L \}, \nonumber
\\
\nonumber & \sum_{l=1}^L u_l' = Y, 
  \\ \nonumber & 0\leq u_l' \leq Y_l', 
 \forall \ l \in \{1,\dots, L \},
   \\ \nonumber & \sum_{l=1}^L \Delta^{S+}_l = \sum_{l=1}^L \Delta^{S-}_l, 
  \\ \nonumber & 0\leq \Delta^{S+}_l \leq \sum_{k\neq l} (s_k)^+, \ \ \  0\leq \Delta^{S-}_l \leq - (s_l)^+, 
  \\\nonumber & \forall \ l \in \{1,\dots, L \}, 
\\ \nonumber & (s_l + \Delta^{S+}_l - \Delta^{S-}_l)\leq y_l \leq (s_l + \Delta^{S+}_l - \Delta^{S-}_l)+U_l +u_l'G, 
\\ & \nonumber \forall \ l \in \{1,\dots, L \}, 
 \\ &\nonumber r_l^n \geq y_l -d_l^n, \ \ o_l^n \geq d_l^n - y_l, \ \forall \ l \in \{1,\dots, L \},
 \\ & \nonumber n \in \{1,\dots, M \}, \text{ and } \nonumber
 \\ & \nonumber r_l^n, \ o_l^n \in \mathbb{Z}^+ \ \forall n \in \{1,\dots, M \}, \ u_l',\ \Delta^{S+}_l, \Delta^{S-}_l, \ y_l \in \mathbb{Z},
 \\& \eta_l, \ \zeta_l \in \mathbb{R}  \ \ \forall \ l \in \{1,\dots, L \}.
 \end{align} 
 This heuristic utilizes significantly fewer integer variables compared to the integer program in \eqref{eq_jr}. Additionally, we have the following result that shows \ac{LSF} can be solved as a linear program to obtain an optimal solution when module capacity equals $1$. This result improves the speed of implementation dramatically in such instances, in comparison with \ac{RSF}.
\begin{proposition} \label{laj_g1} \ac{LSF} can be solved exactly by relaxing all the integrality constraints when module capacity $G=1$. 
\end{proposition}
Proof of this result follows the proof of \autoref{laj_g1} \cite{Malladi2018a}. We remark that the numerical results in \autoref{comput} will justify the robustness of the $G=1$ assumption.
\subsection{\acf{RRO} \label{glr}}
We now consider a distributed decision-making structure in which all the relocation decisions are made using the heuristic \acf{RRO} while replenishment decisions are made at the individual locations. In  \eqref{blend_hr}, consider the inner minimization and note the terms in the inner brackets are bounded below by 
 \begin{align}
& \nonumber \sum_{l=1}^L \bigg [(1-\theta)\min_{y_l} \big\{ \mathcal{L}_l(x,y_l) 
\\ & \nonumber + \beta \sum_{\bm{d},\bm{z}} \sigma(\bm{d},\bm{z},\bm{x}) v^F_l (\bm{\lambda}(\bm{d},\bm{z},\bm{x}), y_l-d_l, Y_l') \big\} 
\\  \nonumber & + \theta \min_{y_l} \big \{ \mathcal{L}_l (\bm{x},y_l) 
\\ \label{eq} &+ \beta \sum_{\bm{d},\bm{z}} \sigma(\bm{d},\bm{z},\bm{x}) v^F_l(\bm{\lambda}(\bm{d},\bm{z},\bm{x}), y_l-d_l, u_l')\big\}  \bigg], \end{align}
where the first minimization is now relaxed to operate over all $y_l$ such that $s_l+\Delta^S_l \leq y_l \leq s_l +\Delta^S_l +U_l  + Y_l' G$ and the second minimization is over all $y_l$ such that $s_l+\Delta^S_l \leq y_l \leq s_l +\Delta^S_l +U_l  + u_l' G$.
 We note that the terms in \eqref{eq} equal
$$  \sum_{l=1}^L \bigg[ (1-\theta) v^F_l(\bm{x},s_l+\Delta^S_l, Y_l') + \theta v^F_l(\bm{x},s_l+\Delta^S_l, u_l') \bigg]. $$

The value function approximation from \autoref{stat_L1} is employed here for tractability of calculations. Then, for the \ac{RRO} heuristic, at every decision epoch with beginning state $(\bm{x}, \bm{s}, \bm{u})$, we first solve
\begin{enumerate}[wide, labelwidth=!, labelindent=0pt]
\item the following integer program to determine the transport decisions, \textit{i.e.}, the amount of inventory $\Delta^S$ received at every location $l$ and the number of production modules $\Delta^M$ received at every location $l$: 
\begin{align}
 & \nonumber \text{\ac{RRO}:} 
 \\  \nonumber & \min \sum_{l=1}^L\sum_{\Delta^S =-(s_l)^+}^{\sum_{k\neq l} (s_k)^+} \sum_{\Delta^M = -u_l }^{Y_l'-u_l}  w(l,\Delta^S, \Delta^M)  \bigg\{ K^{S+}_l  (\Delta^S)^+ +
 \\ &  \nonumber K^{S-}_l (-\Delta^S)^+ + K^M\abs{\Delta^M}/2 
  + {\hat{v}}_l^{F,\theta}(\bm{x}, s_l+\Delta^S, u_l+\Delta^M) \bigg\},
\\ \nonumber & \text{subject to }  
\\ \nonumber &   \sum_{l=1}^L\sum_{\Delta^S =-(s_l)^+}^{\sum_{k\neq l} (s_k)^+} \sum_{\Delta^M = -u_l }^{Y_l'-u_l}  w(l,\Delta^S, \Delta^M) \ \Delta^M = 0, 
 \\ \nonumber &  \sum_{l=1}^L\sum_{\Delta^S =-(s_l)^+}^{\sum_{k\neq l} (s_k)^+} \sum_{\Delta^M = -u_l }^{Y_l'-u_l}  w(l,\Delta^S, \Delta^M) \ \Delta^S =0, \ \text { and }
 \\ \nonumber & \ w(l,\Delta^S, \Delta^M) \in \{0,1\}, \ \forall  \ \Delta^S \in \{ -(s_l)^+, \dots, \sum_{k\neq l} (s_k)^+  \}, 
 \\&  \Delta^M \in \{ -u_l, \dots, Y_l'-u_l \},
 \ l \in \{1,\dots,L\}.
 \end{align} 
\item  We then determine the local controllers' replenishment decisions through the location-wise order-up-to-policy presented in \eqref{eq_ponews}, in which the quantity transshipped to any location $l$ will be obtained using the solution of the above integer program \ac{RRO} as $\Delta^S_l =$  $\sum_{\Delta^S= -(s_l)^+}^{\sum_{k\neq l} (s_k)^+ }$ $ \sum_{\Delta^M = -u_l}^{Y_l'-u_l} w(l,\Delta^S, \Delta^M) \Delta^S$ for all locations $l$.
\end{enumerate}
\section{Computational Study and Results  \label{comput}}
In \autoref{inst_set}, we present the experimental design of generating instances that would allow us to study the variation of heuristic quality and the value added due to mobility as a function of  
the number of modulation states $N$, the probability of not transitioning away from any modulation state $\phi$,
the number of locations $L$, the module capacity $G$, the movement cost per unit of inventory between any pair of locations $K^S$, and the movement cost per production module $K^M$. 
On each instance of the generated instance sets, we implement the heuristic policies proposed in \autoref{heurist} on fifty sampled trajectories to obtain a sample average cost of performance for each policy. All the policies are then compared against the selected benchmark policy, \ac{DNF}. We then present an analysis of our computational findings in  Sections \ref{comput_res} and \ref{valMobility}. 
\subsection{Instance design \label{inst_set}}
We generated two sets of instances in the following manner.
\subsubsection{Set A}
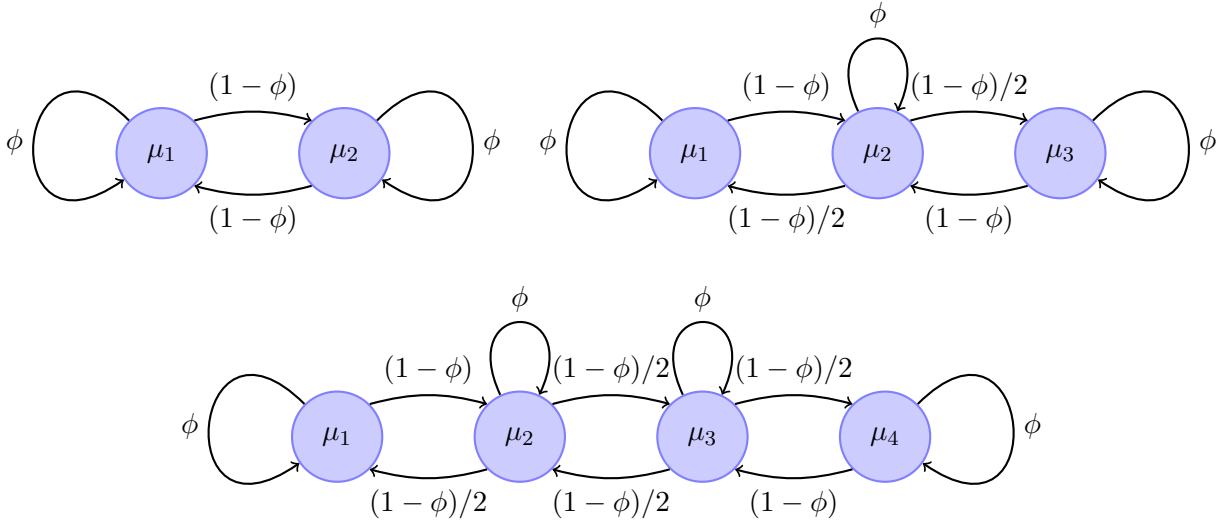
\begin{figure}[h!]
\centering

\begin{tikzpicture}
    [
        node distance =1.2cm,
        place/.style={circle,draw=blue!50,fill=blue!20,thick,
                      inner sep=0pt,minimum size=7mm}
    ]
    \node[place] (1) {$\mu_{1}$};
    \node[place] (2) [right=of 1] {$\mu_{2}$};

\draw [thick,->] (1) edge [looseness=7, out= 135, in=215]node[left]{$\phi$} (1);
\draw [->,thick] (1.north east) to [bend left=15]  node[above] {$(1-\phi)$}  (2.north west);
\draw [<-,thick] (1.south east) to [bend right=15] node[below] {$(1-\phi)$}  (2.south west);
\draw [thick,->] (2) edge [looseness=7, in=325, out = 45] node[right]{$\phi$} (2);
\end{tikzpicture}
\begin{tikzpicture}
    [
        node distance =1.2cm,
        place/.style={circle,draw=blue!50,fill=blue!20,thick,
                      inner sep=0pt,minimum size=7mm}
    ]
    \node[place] (1) {$\mu_{1}$};
    \node[place] (2) [right=of 1] {$\mu_{2}$};
    \node[place] (3) [right=of 2] {$\mu_{3}$};

\draw [thick,->] (1) edge [looseness=7, out= 135, in=215]node[left]{$\phi$} (1);
\draw [->,thick] (2.north east) to [bend left=15]  node[above] {$(1-\phi)/2$}  (3.north west);
\draw [->,thick] (1.north east) to [bend left=15]  node[above] {$(1-\phi)$}  (2.north west);
  \draw [thick,->] (2) edge [looseness=7, in=65, out = 115] node[above]{$\phi$} (2);
\draw [<-,thick] (1.south east) to [bend right=15]  node[below] {$(1-\phi)/2$}  (2.south west);
\draw [<-,thick] (2.south east) to [bend right=15] node[below] {$(1-\phi)$}  (3.south west);
\draw [thick,->] (3) edge [looseness=7, in=325, out = 45] node[right]{$\phi$} (3);
\end{tikzpicture}
\begin{tikzpicture}
    [
        node distance =1.2cm,
        place/.style={circle,draw=blue!50,fill=blue!20,thick,
                      inner sep=0pt,minimum size=7mm}
    ]
    \node[place] (1) {$\mu_{1}$};
    \node[place] (2) [right=of 1] {$\mu_{2}$};
    \node[place] (3) [right=of 2] {$\mu_{3}$};
    \node[place] (4) [right=of 3] {$\mu_{4}$};

\draw [thick,->] (1) edge [looseness=7, out= 135, in=215]node[left]{$\phi$} (1);
\draw [->,thick] (2.north east) to [bend left=15]  node[above] {$(1-\phi)/2$}  (3.north west);
\draw [->,thick] (1.north east) to [bend left=15]  node[above] {$(1-\phi)$}  (2.north west);
  \draw [thick,->] (2) edge [looseness=7, in=65, out = 115] node[above]{$\phi$} (2);
\draw [<-,thick] (1.south east) to [bend right=15]  node[below] {$(1-\phi)/2$}  (2.south west);
\draw [->,thick] (3.north east) to [bend left=15]  node[above] {$(1-\phi)/2$}  (4.north west);
\draw [<-,thick] (2.south east) to [bend right=15] node[below] {$(1-\phi)/2$}  (3.south west);
\draw [thick,->] (3) edge [looseness=7, in=65, out = 115] node[above]{$\phi$} (3);
\draw [<-,thick] (3.south east) to [bend right=15] node[below] {$(1-\phi)$}  (4.south west);
\draw [thick,->] (4) edge [looseness=7, in=325, out = 45] node[right]{$\phi$} (4);
\end{tikzpicture}
\caption{Underlying Markov chain of the modulation process in Set A for $N=2,3,$ and 4 respectively}\label{expt_mc_setA} 
\end{figure}
We fix the following parameters: length of the horizon $T = 30$, number of locations $L = 5$, and total number of modules $Y =  \ceil{\frac{4}{3}L}$. We vary the module capacity $G \in \{1,2,5\}$, fixing the number of demand outcomes $M = 2G+1$ (allowing all integer outcomes between $0$ and $2G$) at each location. We consider three different values for the number of modulation states $N \in \{2, 3, 4 \}$. The underlying Markov chain's transition structure is presented in \autoref{expt_mc_setA}. We vary the probability of not leaving any modulation state, which we refer to  as the \textit{staying probability}, $\phi \in \{0.75,0.95 \}$. We randomly obtain a multi-location discrete demand distribution for each combination of the parameters listed so far, with demand outcomes $\{0, \dots, 2G\}$ such that the probabilities are randomly generated ensuring that exactly one of the $N$ expected demands at each location lies in each of following the intervals 
\begin{itemize}
\item $[0,G)$ and $ [G,2G]$ if $N=2$
\item $[0,0.6G)$, $[0.6G, 1.4G)$, and $ [1.4G,2G]$ if $N=3$, and
\item $[0,0.5G)$, $[0.5G, G)$, $ [G,1.5G]$, and $ [1.5G,2G]$ if $N=4$.
\end{itemize}
We fix the backorder cost $b$ to $2$ and the holding cost to $1$.  
We pair each combination of transshipment cost $K^S \in \{0, 1.5, 2,2.5,10000\}$ and module movement cost $K^M \in \{0, 1.5, 2$, $ 2.5, 10000\}$ with the demand instances created above. 
There are a total of $3 \times 3 \times 2 \times 25 = \bm{450}$ randomly generated instances, with $18$ underlying demand instances.

\subsubsection{Set B}
In this instance set, we focus on varying the number of locations $L \in \{2,5,10,15,$ $20, 25\}$ and the movements costs $K^S$ and $K^M \in \{0, 1.5, 2$, $ 2.5, 1000\}$. We fix the remaining parameters as follows: length of the horizon $T = 30$, total number of modules $Y =  \ceil{\frac{4}{3}L}$, module capacity $G = 1$, number of modulation states $N =3$, staying probability $\phi =0.95 $, number of demand outcomes $M =2G+1$, backorder cost $b =2$, and holding cost $h=1$. A multi-location discrete demand distribution for each combination of these parameters are randomly generated ensuring that exactly one of the $N$ expected demands at each location lies in each of the intervals $[0,0.6G)$, $[0.6G, 1.4G)$, and $ [1.4G,2G]$. This procedure results in a total of $6\times 25 = \bm{150}$ instances, with $6$ underlying demand instances. 
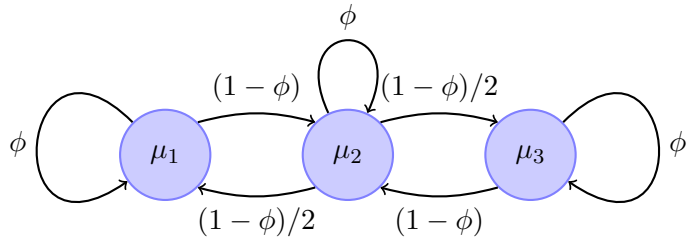
\begin{figure}[h!]
\centering
\begin{tikzpicture}
    [
        node distance =1.2cm,
        place/.style={circle,draw=blue!50,fill=blue!20,thick,
                      inner sep=0pt,minimum size=7mm}
    ]
    \node[place] (1) {$\mu_{1}$};
    \node[place] (2) [right=of 1] {$\mu_{2}$};
    \node[place] (3) [right=of 2] {$\mu_{3}$};

\draw [thick,->] (1) edge [looseness=7, out= 135, in=215]node[left]{$\phi$} (1);
\draw [->,thick] (2.north east) to [bend left=15]  node[above] {$(1-\phi)/2$}  (3.north west);
\draw [->,thick] (1.north east) to [bend left=15]  node[above] {$(1-\phi)$}  (2.north west);
  \draw [thick,->] (2) edge [looseness=7, in=65, out = 115] node[above]{$\phi$} (2);
\draw [<-,thick] (1.south east) to [bend right=15]  node[below] {$(1-\phi)/2$}  (2.south west);
\draw [<-,thick] (2.south east) to [bend right=15] node[below] {$(1-\phi)$}  (3.south west);
\draw [thick,->] (3) edge [looseness=7, in=325, out = 45] node[right]{$\phi$} (3);
\end{tikzpicture}
\caption{Underlying Markov chain of the modulation process in Set B}\label{expt_mc_setB} 
\end{figure}

 Without loss of generality, we set the production cost $c_l$ at all locations to zero in the instances of both sets. 

We note that making the cost parameters location-dependent would make the experiments more realistic. However, since we aim to demonstrate the value of mobile production capacity in interaction with inventory transshipment as a proof of concept even when all the locations are identical, we have generated our instances thus.

\subsection{Quality of Heuristics \label{comput_res}}
We evaluated the heuristic policies, \ac{RRO}, \ac{RSF}, \ac{LSF}, and \ac{MP} on fifty sample trajectories of the instance set, obtained by Monte Carlo simulation, for five values of the coefficient $\theta \in \{0, 0.2, 0.5, 0.8, 1 \}$ where relevant. We compared their performance against the benchmark policy \ac{DNF} and also juxtapose \ac{MNF} against \ac{DNF}.
For each instance, we computed the approximate value function of the $L=1,Y=0$ problem with the various capacities and  determined the minimum total fixed cost among all configurations. We then generated $50$ sample demand trajectories at each epoch based on the current simulated modulation state. For each trajectory, the beginning state is the zero inventory position at all locations and the module configuration that minimizes the sum of the fixed expected total cost of the single location problems with the steady state belief-based distribution of demand as the epoch-invariant demand distribution at each location. We computed the upper bound $ \hat{v}_l(\bm{x},s_l,u_l)$ for $\bm{x} \in X', u \in \{0,\dots, Y\}, \forall \ s_l, \forall \ l$ in a one time offline pre-computation step. We approximated the belief space $X =\{ \bm{x}: \sum_{i=1}^N x_i = 1, x_i \geq 0\}$ with its non-empty, fixed, finite subset $X' = \{ \bm{x}: \sum_{i=1}^N x_i = 1, x_i \in \{0,1/3,2/3,1 \}\} \cup \bm{\pi}$, for $\pi$ such that $\bm{\pi} = \bm{\pi} \hat{P}$ when it exists  \cite{Lovejoy1991}. 

We performed a forward dynamic programming pass or a forward rollout implementing the decision-making proposed by each method at each epoch. We obtained the average performance of each heuristic over the $50$ simulated trajectories of each instance to analyze various resultant trends in comparison to \ac{DNF}.

We compared heuristic performance across values of the coefficient $\theta$ (in Tables \ref{glr-blend-a}, \ref{laj-blend-a}, and \ref{jr-blend-a}  in Appendix \ref{ex_tables}) and found the best performance usually  at $\theta=0.2$ for all the heuristics. \autoref{all_seta_0.2} presents the comparison of the performance of all heuristics at $\theta = 0.2$. We find that although \ac{LSF} at $\theta = 0.2$ is the best performer, for other values of $\theta>0$, \ac{RRO} outperforms it. For $\theta = 0$, the cost of \ac{RRO} is very negative as the inventory holding and backordering components are not considered in the policy's immediate cost.  We note that the cost of the naive policy \ac{MP} is worse than that of \ac{DNF} for $G=1$. However, for higher $G$, \ac{MP} results in significant savings. This observation establishes the need for intelligent, dynamic heuristics that account for future costs, especially when $G=1$.  The proposed heuristics provide about a $38\%-44\%$ average reduction in cost compared to \ac{DNF}, in effect extracting $38\%-44\%$ improvement in system performance from the two forms of mobile flexibility. We note that heuristic quality is almost identical between \ac{RRO}, \ac{LSF}, and \ac{RSF}. 

\begin{table}[h!]
\centering
\caption{Variation of average savings due to proposed heuristics over \ac{DNF} (no flexibility case) with varying $G$ when $\theta = 0.2$ for Instance Set A}
\label{all_seta_0.2}
\begin{tabular}{cccccc}
\toprule
$G $   & \ac{MNF} & \ac{MP} & \ac{RRO} & \ac{LSF} & \ac{RSF}\\
\hline
1       & -9\%  & -5\% & 42\% & 43\% & 44\% \\
2       & -10\% & 13\% & 38\% & 38\% & -  \\
5       & -11\% & 19\% & 44\% & 44\% & -     \\
\hline
Overall & -10\% & 9\%  & 41\% & 42\% &  -  
\\\bottomrule  \hline
\end{tabular}
\end{table}

We repeat the experiments on Set A with a shorter horizon $T=10$ instead of $T=30$ (\autoref{all_seta_0.2_shorter_hor} in Appendix \ref{ex_tables}) and find that \ac{LSF} (which mimics \ac{RSF}) outperforms \ac{RRO} on average by $2\%-3\%$. 
The strength of \ac{RRO} is its unique usefulness while managing instances where different locations are coupled (or correlated) not only through the modulation process. \ac{RSF} and \ac{LSF} rely on the assumption that the demands at different locations are mutually independent, conditional on the belief state.

We now consider the computational efficiency  of the heuristics. \autoref{seta_runtime_G} presents the time taken to compute the policy on a single trajectory using \ac{MP}, \ac{RRO}, \ac{LSF}, and \ac{RSF} and the time taken to compute $\hat{v}^F_l$ for all locations per instance. We note that \ac{MP} is the fastest while \ac{RRO} and \ac{LSF} are significantly faster than \ac{RSF}. Between \ac{RRO} and \ac{LSF}, \ac{RRO} is faster, with a clear edge for $G>1$. 
\begin{table}[h!]
\centering
\caption{Variation of average runtime in seconds with respect to $G$  for Instance Set A}
\label{seta_runtime_G}
\begin{tabular}{cccccc}
\toprule
$G$ & \ac{MP}   & \ac{RRO}  & \ac{LSF}  & \ac{RSF} & $\hat{v}^F_l \ \forall \ l$   \\
\multicolumn{1}{c}{} & \multicolumn{4}{c}{per trajectory}    & per instance        \\
\hline
1 & 0.89 & 1.46 & 2.17 & 2081 & 41   \\
2 & 1.31 & 1.97 & 3.94 &- & 211  \\
5 & 2.61 & 3.64 & 7.59 &- & 2384
\\ \bottomrule \hline
\end{tabular}\end{table}

 \ac{LSF} is computationally faster than \ac{RRO}  on Set B (\autoref{setb_runtime_L}) that contains only $G=1$ instances as \ac{LSF} can be solved as a linear program for $G=1$. 
\begin{table}[h!]
\centering
\caption{Variation of average runtime per trajectory in seconds with respect to $L$  for Instance Set B}
\label{setb_runtime_L}
\begin{tabular}{ccc}
\toprule
$L$       & \ac{RRO}/MP & \ac{LSF}/MP \\
\hline
2       & 0.66   & 2.14   \\
5       & 1.44   & 2.26   \\
10      & 2.97   & 2.05   \\
15      & 6.2    & 1.94   \\
20      & 8.2    & 1.9    \\
25      & 11.8   & 1.7    
\\ \bottomrule \hline    
\end{tabular}
\end{table}

\begin{table}[h!]
\centering
\caption{Variation of average runtime in seconds with respect to $N$  for Instance Set A}
\label{seta_runtime_N}
\begin{tabular}{ccccc}
\toprule
$N$ & \ac{MP}   & \ac{RRO}  & \ac{LSF}  & $\hat{v}^F_l \ \forall \ l$   \\
\multicolumn{1}{c}{} & \multicolumn{3}{c}{per trajectory}    & per instance        \\
\hline
2 & 0.93 & 2.43 & 3.56 & 294  \\
3 & 1.47 & 2.23 & 4.35 & 539  \\
4 & 2.41 & 2.41 & 5.79 & 1802
\\ \bottomrule \hline
\end{tabular}\end{table}
\autoref{seta_runtime_N} shows that the computational effort of computing $\hat{v}^F_l\ \forall \ l$ increases significantly when the number of belief states considered ($N$) is increased.

Thus, we find that \ac{LSF} works best for for $G=1$ and \ac{RRO} for $G>1$ when considering both speed and performance. 

\subsection{Value of Mobility \label{valMobility}}
We now study the trends of value addition due to mobility of production capacity. \autoref{all_seta_0.2} indicates that with increase in $G$, the percentage of savings does not show a clear trend. This behavior might be partly attributed to the fact that the movement cost per unit of capacity is lower when $G$ is higher for the same movement cost per module and the amount of free capacity per location does not scale properly with problem size. 

\begin{table}[h!]
\centering
\caption{Variation of average savings due to \ac{LSF} with $\theta = 0.2$ over \ac{DNF} (no flexibility case) with varying $N$  for Instance Set A}
\label{laj_nmod_seta}
\begin{tabular}{ccccc}
\toprule
$N$       & Overall &  $\phi = 0.95$ &  $\phi = 0.95$ & $\phi = 0.95$ \\
      & &  &  $G=1$ & $G=5$ \\
\hline
2 & 44\% & 44\% & 44\% & 44\% \\
3 & 35\% & 52\% & 28\% & 65\% \\
4 & 46\% & 54\% & 59\% & 56\%      \\
\hline
Overall & 42\%    & 50\%           & 44\%           & 55\%  
\\\bottomrule  \hline      
\end{tabular}
\end{table}

When the number of modulation states $N$ is varied, the average savings over \ac{DNF} due to \ac{LSF} do not exhibit a clear trend but we note that the configuration of the other parameters, such as staying probability $\phi$  and module capacity $G$, affect the influence of $N$ (\autoref{laj_nmod_seta}) on the amount of savings.  It is interesting to note the profit potential in certain configurations: when $G=5$, $N=3$, and $\phi = 0.95$, mobility extracts about 65\% savings. With respect to the probability of not leaving in any modulation state $\phi$, we find that when the dynamics of the world are such that $\phi$ is closer to $1$, about $17\%$ higher average savings are observed (\autoref{all_phi_seta}) than when it is farther.  
\begin{table}[h!]
\centering
\caption{Variation of average savings over \ac{DNF} (no flexibility case) with varying $\phi$  for Instance Set A}
\label{all_phi_seta}
\begin{tabular}{cccc}
\toprule
$\phi$  & \ac{MP}   & \ac{RRO} ($\theta = 0.2$)  & \ac{LSF} ($\theta = 0.2$)  \\
\hline
0.75 & 1\%  & 33\% & 33\% \\
0.95 & 17\% & 48\% & 50\%
\\\bottomrule  \hline  
\end{tabular}
\end{table}

\autoref{val_laj_seta} presents the value of mobility expressed as percentage savings due to \ac{LSF} over \ac{DNF} as a function of movement costs $K^S$ and $K^M$. We note that both forms of flexibility offer significantly high savings even when operated independently as seen from the row / column with the cost set to 1,000. We note that for all the considered combinations of movement costs (except 1,000 for both), significantly high savings, to the tune of 40\%, are observed. We also make note that production capacity mobility independently  ($K^S$ = 1,000) extracts $3-5\%$ higher savings than transshipment operated independently ($K^M$ = 1,000), emphasizing the value of production capacity mobility in comparison to transshipment.  For the subset of Set A with $G=5$, we note that the independent savings from production capacity mobility are about 10-13\% higher than those from inventory mobility; these quantities are almost twice those at $G=1$ (\autoref{val_laj_seta_g1} and \autoref{val_laj_seta_g5} in Appendix \ref{ex_tables}).  
\begin{table}[h!]
\centering
\caption{Value of mobility represented as \% savings using the heuristic \ac{LSF} over \ac{DNF} (no flexibility case) with $\theta = 0.2$ across varying $ K^S$ and $K^M$  for Instance Set A}
\label{val_laj_seta}
\begin{tabular}{cc|ccccc}
\toprule
                         &     & \multicolumn{5}{c}{Module movement cost $K^M$}  \\
                           &  &0   & 1.5      & 2        & 2.5 &1000       \\ 
\midrule
&0    & 49\% & 48\% & 46\% & 50\% & 50\% \\
&1.5  & 48\% & 41\% & 40\% & 42\% & 41\% \\
&2    & 49\% & 41\% & 41\% & 41\% & 38\% \\
&2.5  & 47\% & 43\% & 41\% & 39\% & 37\% \\
\multirow{-6}{*}{\rotatebox[origin=c]{90}{\parbox{2.3cm}{\centering \scriptsize{Transshipment} \\  cost $K^S$}}} &1000 & 50\% & 44\% & 43\% & 40\% & -4\%
\\ \bottomrule  \hline
\end{tabular}
\end{table}

From \autoref{all_L_setb} that presents a trend of average savings from \ac{LSF} in Set B containing only $G=1$ instances, we note that as the number of locations $L$ increases, the average value addition due to resource mobility over \ac{DNF} is very high (30-65\%) generally. Once again, certain configurations extract very high savings. Although an increasing trend is expected as seen in \cite{Malladi2018a}, we do not see it clearly in these averages. 
\begin{table}[h!]
\centering
\caption{ Variation of average savings due to \ac{LSF} with $\theta = 0.2$ over \ac{DNF} (no flexibility case) with varying $L$  for Instance Set B}
\label{all_L_setb}
\begin{tabular}{cccccccc}
\toprule
$L$ &  2    & 5    & 10   & 15   & 20   & 25   & Overall \\
\hline
LAJ & 57\% & 44\% & 32\% & 61\% & 64\% & 48\% & 51\% 
\\\bottomrule  \hline  
\end{tabular}
\end{table}

Since the efficiency of the heuristics decreases with increase in the number of modulation states $N$ in partially observed (PO) decision-making, we consider the case where the decision maker uses the steady state (SS) distribution $\pi$ (when it exists) as the \textit{epoch-invariant} belief state without dynamics. From \autoref{po_ss_co_seta}, we note that \ac{LSF} with SS and \ac{LSF} with PO decision-making perform almost identically. However, on shorter horizons, the additional savings from PO over SS are about 2-3\% higher (\autoref{all_seta_0.2_shorter_hor} of \autoref{ex_tables}). This observation reinforces the value addition of non-stationary demand modeling over short horizons although over long horizons, accurate stationary demand distributions would yield similar savings.  We also compare PO with the case where the modulation process is completely observed and find that complete observability of the modulation process improves overall savings by only 1\%.   

\begin{table}[h!]
\centering
\caption{Comparison of average savings \ac{LSF} with $\theta = 0.2$ over \ac{DNF} (no flexibility case) when a) the DM models epoch-invariant steady state (SS) demand distributions when a partially observed (PO) modulation process is acting, b) the DM partially observes the modulation process, and c) the DM completely observes (CO) the modulation process  for Instance Set A}
\label{po_ss_co_seta}
\begin{tabular}{cccc}
\toprule
$G$       & SS   & PO   & CO   \\
\hline
1       & 44\% & 43\% & 44\% \\
2       & 38\% & 38\% & 39\% \\
5       & 44\% & 44\% & 44\% \\
\hline
Overall & 42\% & 42\% & 42\%
\\\bottomrule  \hline
\end{tabular}
\end{table}

Thus, we conclude our computational analysis by emphasizing the significant impact of module mobility, the value of module capacity, and nature of epoch-variance of demands on the performance of production-inventory systems under epoch-variant demands and while indicating the comparative advantages of the heuristics \ac{LSF} and \ac{RRO} over \ac{RSF}. In particular, \ac{LSF} for $G=1$ and \ac{RRO} for $G>1$ are efficient as well as effective. 

\section{Conclusion \label{concl}}
We have modeled a multi-location production-inventory system under stochastic demand with geographically relocatable production capacity and have developed computationally efficient heuristic methods for this problem.  We have shown that the heuristics \ac{LSF} and \ac{RRO} are computationally efficient and improve in solution quality as the system size and uncertainty increase for the $L$ location, $Y$ module problem. We have observed the value of mobility of production capacity to be around 41\% on average, relative to systems with no production capacity mobility, irrespective of the presence of transshipment flexibility. We have noted that making decisions assuming a stationary belief state set to the steady state distribution $\pi$ performs comparably with decision-making with partially observed Markov-modulated demands. Complete observability of the modulation state does not appear to add significant value in the current context. Additionally, we infer that although centralized control (\ac{RSF} and \ac{LSF}) results in slightly lower costs, decentralized control (\ac{RRO}) heuristics perform significantly faster for $G>1$.

We conclude that data-driven production capacity relocation represents a promising new supply chain design and operations feature, and we have presented effective heuristics to determine inventory replenishment and production capacity relocation, with or without inventory transshipment. Our results serve to justify more detailed analyses for specific cases having more realistic assumptions regarding frequency of decision epochs, lead times, and transshipment costs and that this analysis could use the heuristics developed in this paper for giving guidance to simulation models.

\bibliographystyle{spbasic-no-url}      
\bibliography{Mendeley_MMPTT}

\begin{thebibliography}{43}
\providecommand{\natexlab}[1]{#1}
\providecommand{\url}[1]{{#1}}
\providecommand{\urlprefix}{URL }
\expandafter\ifx\csname urlstyle\endcsname\relax
  \providecommand{\doi}[1]{DOI~\discretionary{}{}{}#1}\else
  \providecommand{\doi}{DOI~\discretionary{}{}{}\begingroup
  \urlstyle{rm}\Url}\fi
\providecommand{\eprint}[2][]{\url{#2}}

\bibitem[{Axs{\"{a}}ter et~al.(2002)Axs{\"{a}}ter, Marklund, and
  Silver}]{Axsater2002}
Axs{\"{a}}ter S, Marklund J, Silver EA (2002) {Heuristic methods for
  centralized control of one-warehouse , N -retailer inventory systems}.
  Manufactuirng {\&} Service Operations Management 4(1):75--97,
  \doi{10.1287/msom.4.1.75.291}

\bibitem[{{Bayer Technology Services
  GMBH}(2014)}]{BayerTechnologyServicesGMBH2014}
{Bayer Technology Services GMBH} (2014) {Flexible, fast and future production
  processes}. Tech. rep.

\bibitem[{Bernstein and DeCroix(2006)}]{Bernstein2006}
Bernstein F, DeCroix GA (2006) {Inventory policies in a decentralized assembly
  system}. Operations Research 54(2):324--336, \doi{10.1287/opre.1050.0256}

\bibitem[{Bernstein and Federgruen(2005)}]{Bernstein2005}
Bernstein F, Federgruen A (2005) {Decentralized supply chains with competing
  retailers rnder demand uncertainty}. Management Science 51(1):18--29,
  \doi{10.1287/mnsc.1040.0218}

\bibitem[{Bertsekas et~al.(1997)Bertsekas, Tsitsiklis, and Wu}]{bertsekas97}
Bertsekas DP, Tsitsiklis JN, Wu C (1997) {Rollout algorithms for combinatorial
  optimization}. Journal of Heuristics 3:245--262

\bibitem[{Burnetas and Katehakis(1997)}]{Apostolos1997}
Burnetas AN, Katehakis MN (1997) Optimal adaptive policies for {M}arkov
  decision processes. Mathematics of Operations Research 22(1):222--255,
  \doi{10.1287/moor.22.1.222}

\bibitem[{Cheung and Simchi-Levi(2019)}]{cheung15}
Cheung WC, Simchi-Levi D (2019) Sampling-based approximation schemes for
  capacitated stochastic inventory control models. Mathematics of Operations
  Research 44(2):668--692, \doi{10.1287/moor.2018.0940}

\bibitem[{Federgruen and Zipkin(1986)}]{Federgruen1986}
Federgruen A, Zipkin P (1986) {An inventory model with limited production
  capacity and uncertain demands II: The discounted-cost criterion}.
  Mathematics of Operations Research 11(2):208--215, \doi{10.2307/3689804}

\bibitem[{{Geek Wire}(2018)}]{GeekWire2018}
{Geek Wire} (2018) {Amazon finally wins a patent for 3-D printing on demand,
  for pickup or delivery}

\bibitem[{Ghiani et~al.(2002)Ghiani, Guerriero, and Musmanno}]{ghiani2002}
Ghiani G, Guerriero F, Musmanno R (2002) {The capacitated plant location
  problem with multiple facilities in the same site}. Computers {\&} Operations
  Research 29(13):1903--1912

\bibitem[{Godfrey and Powell(2001)}]{godfrey01}
Godfrey GA, Powell WB (2001) {An adaptive, distribution-free algorithm for the
  newsvendor problem with censored demands, with applications to inventory and
  distribution}. Management Science 47(8):1101--1112

\bibitem[{Goodson et~al.(2017)Goodson, Thomas, and Ohlmann}]{Goodson2017}
Goodson JC, Thomas BW, Ohlmann JW (2017) {A rollout algorithm framework for
  heuristic solutions to finite-horizon stochastic dynamic programs}. European
  Journal of Operational Research 258(1):216--229,
  \doi{10.1016/j.ejor.2016.09.040}

\bibitem[{Halper and Raghavan(2011)}]{Halper2011}
Halper R, Raghavan S (2011) {The mobile facility routing problem}.
  Transportation Science 45(3):413--434, \doi{10.1287/trsc.1100.0335}

\bibitem[{Herer et~al.(2002)Herer, Tzur, and Y{\"{u}}cesan}]{Herer2002}
Herer YT, Tzur M, Y{\"{u}}cesan E (2002) {Transshipments: An emerging inventory
  recourse to achieve supply chain leagility}. International Journal of
  Production Economics 80(3):201--212, \doi{10.1016/S0925-5273(02)00254-2}

\bibitem[{Herer et~al.(2006)Herer, Tzur, and Y{\"{u}}cesan}]{Herer2006}
Herer YT, Tzur M, Y{\"{u}}cesan E (2006) {The multilocation transshipment
  problem}. IIE Transactions (Institute of Industrial Engineers)
  38(3):185--200, \doi{10.1080/07408170500434539}

\bibitem[{Jena et~al.(2015)Jena, Cordeau, and Gendron}]{jena2015}
Jena SD, Cordeau JF, Gendron B (2015) {Dynamic facility location with
  generalized modular capacities}. Transportation Science 49(3):489--499

\bibitem[{Jordan and Graves(1995)}]{Jordan1995}
Jordan WC, Graves SC (1995) {Principles on the benefits of manufacturing
  process flexibility}. Management Science 41(4):577--594,
  \doi{10.1287/mnsc.41.4.577}

\bibitem[{Karmarkar(1979)}]{Karmarkar1979}
Karmarkar US (1979) {Convex/Stochastic programming and multilocation inventory
  problems}. Naval Research Logistics 26(1):1--19, \doi{10.1002/nav.3800260102}

\bibitem[{Karmarkar(1981)}]{Karmarkar1981}
Karmarkar US (1981) {The multiperiod multilocation inventory problem}.
  Operations Research 29(2):215--228

\bibitem[{Karmarkar(1987)}]{Karmarkar1987}
Karmarkar US (1987) {The multilocation multiperiod inventory problem: Bounds
  and approximations}. Management Science 33(1):86--94,
  \doi{10.1287/mnsc.33.1.86}

\bibitem[{Katehakis et~al.(2015)Katehakis, Melamed, and Shi}]{Katehakis2015}
Katehakis MN, Melamed B, Shi JJ (2015) Optimal replenishment rate for inventory
  systems with compound poisson demands and lost sales: a direct treatment of
  time-average cost. Annals of Operations Research
  \doi{10.1007/s10479-015-1998-y}

\bibitem[{Kouvelis and Gutierrez(1997)}]{Kouvelis1997TheSystem}
Kouvelis P, Gutierrez GJ (1997) {The newsvendor problem in a global market:
  Optimal centralized and decentralized control policies for a two-market
  stochastic inventory system}. Management Science 43(5):571--585,
  \doi{10.1287/mnsc.43.5.571}

\bibitem[{Lien et~al.(2011)Lien, Iravani, Smilowitz, and Tzur}]{Lien2011}
Lien RW, Iravani SM, Smilowitz K, Tzur M (2011) {An efficient and robust design
  for transshipment networks}. Production and Operations Management
  20(5):699--713, \doi{10.1111/j.1937-5956.2010.01198.x}

\bibitem[{Lovejoy(1991)}]{Lovejoy1991}
Lovejoy WS (1991) {Computationally feasible bounds for partially observed
  Markov decision processes}. Operations Research 39(1):162--175

\bibitem[{Malladi et~al.(2018)Malladi, Erera, and White~III}]{Malladi2018}
Malladi SS, Erera AL, White~III CC (2018) {Inventory control with modulated
  demand and a partially observed modulation process}. arXiv

\bibitem[{Malladi et~al.(2020)Malladi, Erera, and White~III}]{Malladi2018a}
Malladi SS, Erera AL, White~III CC (2020) A dynamic mobile production capacity
  and inventory control problem. IISE Transactions 52(8):926--943,
  \doi{10.1080/24725854.2019.1693709}

\bibitem[{Marcotte and Montreuil(2016)}]{Marcotte2016IntroducingProduction}
Marcotte S, Montreuil B (2016) {Introducing the concept of hyperconnected
  mobile production}. Progress in Material Handling Research

\bibitem[{Melo et~al.(2005)Melo, Nickel, and da~Gama}]{melo2006}
Melo MT, Nickel S, da~Gama F (2005) {Dynamic multi-commodity capacitated
  facility location: a mathematical modeling framework for strategic supply
  chain planning}. Computers {\&} Operations Research 33(1):181--208

\bibitem[{{MIT News}(2016)}]{MIT2016}
{MIT News} (2016) {Pharmacy on demand}

\bibitem[{{Pfizer}(2015)}]{Pfizer2015}
{Pfizer} (2015) {Pfizer announces collaboration with GSK on next-generation
  design of portable, continuous, miniature and modular (PCMM) oral solid dose
  development and manufacturing units}

\bibitem[{Powell(2007)}]{adpbook}
Powell WB (2007) Approximate dynamic programming: Solving the curses of
  dimensionality (Wiley Series in Probability and Statistics), 2nd edn. John
  Wiley {\&} Sons, Inc., Hoboken, New Jersey, \doi{10.1002/9781118029176}

\bibitem[{Powell(2012)}]{Powell2012}
Powell WB (2012) {Perspectives of approximate dynamic programming}. Annals of
  Operations Research pp 1--38, \doi{10.1007/s10479-012-1077-6}

\bibitem[{Puterman(1994)}]{Puterman1994}
Puterman ML (1994) {Markov decision processes: Discrete stochastic dynamic
  programming}. \doi{10.1080/00401706.1995.10484354}

\bibitem[{Qiu and Sharkey(2013)}]{Qiu2013}
Qiu J, Sharkey TC (2013) {Integrated dynamic single-facility location and
  inventory planning problems}. IIE Transactions (Institute of Industrial
  Engineers) 45(8):883--895, \doi{10.1080/0740817X.2013.770184}

\bibitem[{Rudi et~al.(2001)Rudi, Kapur, and Pyke}]{Rudi2001AMaking}
Rudi N, Kapur S, Pyke DF (2001) {A Two-location inventory model with
  transshipment and local decision making}. Management Science
  47(12):1668--1680, \doi{10.1287/mnsc.47.12.1668.10235}

\bibitem[{Ryzhov et~al.(2012)Ryzhov, Powell, and Frazier}]{Ryzhov2012}
Ryzhov IO, Powell WB, Frazier PI (2012) {The knowledge gradient algorithm for a
  general class of online learning problems}. Operations Research
  60(1):180--195, \doi{10.1287/opre.1110.0999}

\bibitem[{Secomandi(2001)}]{Secomandi2001}
Secomandi N (2001) {A Rollout policy for the vehicle routing problem with
  stochastic demands}. Operations Research 49(5):796--802,
  \doi{10.1287/opre.49.5.796.10608}

\bibitem[{Smallwood and Sondik(1973)}]{Smallwood1973}
Smallwood RD, Sondik EJ (1973) {The optimal control of partially observable
  Markov processes over a finite horizon}. Operations Research
  21(5):1071--1088, \doi{10.1287/opre.21.5.1071}

\bibitem[{Sondik(1978)}]{Sondik1978TheCosts}
Sondik EJ (1978) {The Optimal control of partially observable Markov processes
  over the infinite horizon: discounted costs}. Operations Research
  26(2):282--304, \doi{10.1287/opre.26.2.282}

\bibitem[{Verlinde et~al.(2014)Verlinde, Macharis, Milan, and
  Kin}]{Verlinde2014}
Verlinde S, Macharis C, Milan L, Kin B (2014) {Does a Mobile Depot Make Urban
  Deliveries Faster, More Sustainable and More Economically Viable: Results of
  a Pilot Test in Brussels}. In: Transportation Research Procedia: Mobil. TUM
  2014 “Sustainable Mobility in Metropolitan Regions”, May 19-20, 2014,
  \doi{10.1016/j.trpro.2014.11.027}

\bibitem[{Wee and Dada(2005)}]{Wee2005}
Wee KE, Dada M (2005) {Optimal policies for transshipping inventory in a retail
  network}. Management Science 51(10):1519--1533, \doi{10.1287/mnsc.1050.0441}

\bibitem[{W{\"{o}}rsd{\"{o}}rfer and Lier(2017)}]{worsdorfer2017b}
W{\"{o}}rsd{\"{o}}rfer D, Lier S (2017) {Optimized modular production networks
  in the process industry}. In: Operations Research Proceedings 2015, Springer
  International Publishing

\bibitem[{W{\"{o}}rsd{\"{o}}rfer et~al.(2017)W{\"{o}}rsd{\"{o}}rfer, Lier, and
  Crasselt}]{Worsdorfer2017}
W{\"{o}}rsd{\"{o}}rfer D, Lier S, Crasselt N (2017) {Real options-based
  evaluation model for transformable plant designs in the process industry}.
  Journal of Manufacturing Systems 42:29--43, \doi{10.1016/j.jmsy.2016.11.001}

\end{thebibliography}
\newpage
\section*{Appendices}
\ref{struc_L1} provides the foundational results for the $L = 1$ case on which bounds presented in \autoref{approx_single} for the general $(L, Y)$ case are based. \ref{proof_prop7} presents a proof of Proposition 1, \ref{ex_hr} presents a heuristic that is analogous to the heuristic \ac{LSF}, and \ref{ex_tables} presents additional tables of computational results. 

\renewcommand{\thesubsection}{A\arabic{subsection}} 
\subsection{{Analysis for the $L=1$ Case} \label{struc_L1}}
Assume $v_0 = 0$, $v_{n+1}  = Hv_n$, define $\mathcal{G}_n(\bm{x},y)= \mathcal{G}(\bm{x}, y, v_n)$ for all $n$, and let $y_n^*(\bm{x}, C)$ be the smallest value that minimizes $\mathcal{G}_n(\bm{x}, y)$ with respect to $y$. We remark that
\begin{eqnarray*}
v_{n+1}(\bm{x},s,C) = \begin{cases} \mathcal{G}_n(\bm{x}, s) & \text{ if } s\geq y_n^*(\bm{x},C) 
\\ \mathcal{G}_n(\bm{x},s+C) & \text{ if } s\leq y_n^*(\bm{x},C) - C
\\ \mathcal{G}_n(\bm{x},y_n^*(\bm{x},C)) & \text{ otherwise. } \end{cases}
\end{eqnarray*}
We now present claims for  structured results with respect to $\mathcal{G}_n$, $v_n$, and $y_n^*$ based on results in \cite{Federgruen1986} and \cite{Malladi2018}. 
\begin{proposition}\label{prop1}
For all $n$, $\bm{x}$, and $C$, 
\begin{enumerate}[label = (\roman*)] 
\item \label{item1i} $\mathcal{G}_n(\bm{x},y)$ is convex in $y$
\item \label{item1ii} $v_n(\bm{x},s,C)$ is:
\begin{enumerate}[label = (\alph*)]
\item convex in $s$,
\item non-decreasing for $s\geq y_n^*(\bm{x},C)$,
\item non-increasing for $s\leq y_n^*(\bm{x},C) - C$,
\item equal to $v_n(\bm{x},y_n^*(\bm{x},C), C)$ otherwise
\end{enumerate}
\item \label{item1iii} $v_{n+1}(\bm{x},s,C) \geq v_n(\bm{x},s,C)$ for all $s$. 
\end{enumerate}
\end{proposition}
\proof [\textbf{Proof of Proposition \ref{prop1}}]
The convexity of $\mathcal{G}_0(\bm{x},y)$ in $y$ for all $\bm{x}$ follows from the definitions and assumptions. Assume $\mathcal{G}_n(\bm{x},y)$ is convex in $y$ for all $\bm{x}$. It is then straightforward to show that \autoref{item1ii} holds for $n=n+1$ and all $(\bm{x},C)$. We remark that the function $g(y) = w(f(y))$ is convex and non-decreasing (non-increasing) if $w$ is convex and non-decreasing (non-increasing) and if $f$ is linear and non-decreasing. Hence, $\mathcal{G}_{n+1}(\bm{x},y)$ is convex in $y$ for all $\bm{x}$, and \autoref{item1i} and \autoref{item1ii} hold for all $n$ by induction. Since $v_1(\bm{x},s,C) \geq v_0(\bm{x},s,C)$, a standard induction argument guarantees that \autoref{item1iii} holds.
\endproof
Let $v_n(\bm{x},s)  = v_n(\bm{x},s,C)$, $v_n'(\bm{x},s) = v_n(\bm{x},s,C')$, $\mathcal{G}_n(\bm{x},y) = \mathcal{G}(\bm{x},y,v_n)$, and $\mathcal{G}_n'(\bm{x},y) = \mathcal{G}(\bm{x},y,v_n')$. 
\begin{proposition} \label{prop3}
Assume $C\leq C'$, and that $y_n^*(\bm{x},C) - d  \leq y_n^*(\bm{\lambda}(\bm{d},\bm{z},\bm{x}),C)$ for all $n$ and all $(\bm{d},\bm{z},\bm{x})$. Then for all $n$, $\bm{x}$, and $s$, 
\begin{enumerate}[label = (\roman*)]
\item \label{item3i} $v_n'(\bm{x},s,C) \leq v_n(\bm{x},s,C)$
\item \label{item3ii}If  $y\leq y' \leq y_n^*(\bm{x},C)$, then 
$\mathcal{G}_n(\bm{x},y') - \mathcal{G}_n(\bm{x},y) \leq \mathcal{G}_n'(\bm{x},y') - \mathcal{G}_n'(\bm{x},y)$
\item \label{item3iii}If $s\leq s' \leq y_n^*(\bm{x},C)$, then 
$ v_{n+1}(\bm{x},s',C) - v_{n+1}(\bm{x},s) \leq v_{n+1}'(\bm{x},s',C) - v_{n+1}'(\bm{x},s,C). $
\item \label{item3iv}$y_n^*(\bm{x},C') \leq y_n^*(\bm{x},C)$. 
\end{enumerate}
\end{proposition}
\proof[\textbf{Proof of Proposition \ref{prop3}}]
Proof of \autoref{item3i} is straightforward. Regarding \autoref{item3ii}-\autoref{item3iv}, note \autoref{item3ii} holds for $n=0$; assume \autoref{item3ii} holds for $n$. Then \autoref{item3iv} also holds for $n$. We now outline the proof that \autoref{item3iii} holds for $n=n+1$. Recall
$$ v_{n+1}(\bm{x},s,C) = \begin{cases} \mathcal{G}_n(\bm{x}, s+C)  & \text{ if } s\leq y_n - C 
\\  \mathcal{G}_n(\bm{x},s) & \text{ if } s\geq y_n
\\ \mathcal{G}_n(\bm{x},y_n) & \text{ otherwise,}\end{cases}$$
where $y_n = y_n^*(\bm{x},C)$, and 
$$ v_{n+1}'(\bm{x},s,C) = \begin{cases} \mathcal{G}_n'(\bm{x}, s+C')  & \text{ if } s\leq y_n' - C'
\\  \mathcal{G}_n'(\bm{x},s) & \text{ if } s\geq y_n'
\\ \mathcal{G}_n'(\bm{x},y_n') & \text{ otherwise,}\end{cases}$$
where $y_n' = y_n^*(\bm{x},C')$. Similar to the proof of \autoref{prop2} and the proof of \cite[Theorem ~3]{Federgruen1986}, there are two cases: (1) $y_n - C \leq y_n'$, (2) $y_n' \leq y_n - C$, which are more completely described as 
$$y_n' - C' \leq y_n - C \leq y_n' \leq y_n, $$
$$y_n' - C' \leq y_n' \leq y_n- C \leq y_n,$$
respectively. For each case, there are 10 different sets of inequalities that the pair $(s,s')$ can satisfy. Showing that \autoref{item3iii} holds when $n=n+1$ for each of the 20 sets of inequalities is tedious but straightforward. 
We now show that for $s \leq s'$, 
$$ v_{n+1}(\bm{x},s',C)-v_{n+1}(\bm{x},s,C) \leq v_{n+1}'(\bm{x},s',C) - v_{n+1}'(\bm{x},s,C)$$ 
implies that for $y\leq y' \leq y_n$, $\mathcal{G}_{n+1}(\bm{x},y') - \mathcal{G}_{n+1}(\bm{x},y) \leq \mathcal{G}_{n+1}'(\bm{x},y') -\mathcal{G}_{n+1}'(\bm{x},y)$. Note 
\begin{multline*}v_{n+1}(\bm{\lambda}(d,z,\bm{x}),y'-d,C) - v_{n+1}(\bm{\lambda}(d,z,\bm{x}),y-d,C) 
\\ \leq v_{n+1}'(\bm{\lambda}(d,z,\bm{x}),y'-d,C) - v_{n+1}'(\bm{\lambda}(d,z,\bm{x}),y-d,C) \end{multline*}
for $y-d \leq y'-d \leq y_n^*(\bm{\lambda}(d,z,\bm{x}),C)$, which implies 
$$ \mathcal{G}_{n+1}(\bm{x},y') - \mathcal{G}_{n+1}(\bm{x},y) \leq  \mathcal{G}_{n+1}'(\bm{x},y') - \mathcal{G}_{n+1}'(\bm{x},y) $$
for all $y \leq y' \leq y_{n+1}^*(\bm{x},C)$ assuming $y_{n+1}^*(\bm{x},C) - d \leq y_{n+1}^*(\bm{\lambda}(d,z,\bm{x}), C)$ for all $(d,z,\bm{x})$. A standard induction argument completes the proof. 
\endproof

\begin{proposition} \label{prop2} Assume $y_n^*(\bm{x},C)-d_l \leq y_n^*(\bm{\lambda}(\bm{d},\bm{z},\bm{x}), C)$ for all $n$ and all $(\bm{d},\bm{z},\bm{x})$. Then for all $n$, $s\leq s' \leq y_n^*(\bm{x},C)$ implies:
\begin{enumerate}[label = (\roman*)]
\item \label{item2i}$v_n(\bm{x},s',C) - v_n(\bm{x},s,C) \geq v_{n+1}(\bm{x},s',C) - v_{n+1}(\bm{x},s,C)$,
\item \label{item2ii} $\mathcal{G}_n(\bm{x},s')-\mathcal{G}_n(\bm{x},s) \geq \mathcal{G}_{n+1}(\bm{x},s') - \mathcal{G}_{n+1}(\bm{x},s)$,
\item \label{item2iii}$y_n^*(\bm{x},C) \leq y_{n+1}^*(\bm{x},C)$. 
\end{enumerate} 
\end{proposition}

\proof[\textbf{Proof of Proposition \ref{prop2}}]
We note \autoref{item2i} holds when $n=0$. Assume \autoref{item2i} holds for $n=n-1$. Let $y\leq y' \leq y_{n-1}^*(\bm{x}, C)$, implying that $y-d \leq y'-d \leq y_{n-1}^*(\bm{x},C)-d  \leq y^*_{n-1}(\bm{\lambda}(d,z,\bm{x}),C)$ for all $(d,z,\bm{x})$. Hence, 
\begin{multline*}v_{n-1} (\bm{\lambda}(d,z,\bm{x}),y'-d,C) - v_{n-1}(\bm{\lambda}(d,z,\bm{x}),y-d, C) 
\\ \geq v_n (\bm{\lambda}(d,z,\bm{x}),y'-d, C) - v_n(\bm{\lambda}(d,z,\bm{x}),y-d, C), \end{multline*}
and thus \autoref{item2ii} holds for $n=n-1$ for all $y\leq y' \leq y_{n-1}^*(\bm{x},C)$. Letting $y' = y_{n-1}^*(\bm{x},C)$, we observe 
\begin{multline*}0 \geq \mathcal{G}_{n-1}(\bm{x},y^*_{n-1}(\bm{x},C)) - \mathcal{G}_{n-1}(\bm{x},y) 
\\ \geq \mathcal{G}_n(\bm{x}, y_{n-1}^*(\bm{x},C)) - \mathcal{G}_n(\bm{x},y);\end{multline*}
hence, \autoref{item2iii} holds for $n=n-1$. 

We now outline a proof that $s\leq s' \leq y_n^*(\bm{x},C)$ implies 
\begin{eqnarray}\label{eqstar} v_n(\bm{x},s') - v_n(\bm{x},s) \geq v_{n+1}(\bm{x},s') - v_{n+1}(\bm{x},s).\end{eqnarray}
Following an argument in the proof of \cite[Theorem~2]{Federgruen1986}, we consider two general cases: (1) $y_n^*(\bm{x},C)-C \leq y_{n-1}^*(\bm{x},C)$ and (2) $y_{n-1}^*(\bm{x},C) \leq y_n^*(\bm{x},C)-C$. Letting the dependence on $(\bm{x},C)$ be implicit, cases (1) and (2) are more completely described as 
$$ y_{n-1}^*- C \leq y_n^* - C \leq y_{n-1}^* \leq y_n^* $$
$$ y_{n-1}^*- C \leq y_{n-1}^* \leq y_n^* - C \leq y_n^*, $$
respectively. For each case, there are $10$ different sets of inequalities that the pair $(s,s')$ can satisfy. The values $v_n(\bm{x},s'), v_n(\bm{x},s), v_{n+1}(\bm{x},s')$, and $v_{n+1}(\bm{x},s)$ are well defined for each of these inequalities in terms of $\mathcal{G}_{n-1}$ and $\mathcal{G}_n$. Showing that \eqref{eqstar} holds for each of these 20 different sets of inequalities is again tedious but straightforward. 

A standard induction argument completes the proof of the proposition.
\endproof

We now claim that $v(\bm{x},s,C)$ is convex in $C$. 
\begin{proposition} \label{prop4}
\begin{enumerate}[label = (\roman*)]
\item \label{item4i}If $y\in A(s,C)$ and $y' \in A(s, C')$, then $\lambda y + (1-\lambda)y'$ $\in A(s, \lambda C+(1-\lambda)C')$. 
\item \label{item4ii}If $\xi \in A(s, \lambda C+(1-\lambda)C')$, then there is a $y \in A(s,C)$ and a $y' \in A(s,C')$ such that $\xi =\lambda y+ (1-\lambda)y'$.
\item \label{item4iii}For real-valued and continuous $v$, 
\begin{multline*}
\min\{v(\xi): \xi \in A(s, \lambda C+(1-\lambda)C') \}
\\  = \min\{ v(\lambda y+(1-\lambda)y'):y\in A(s,C) \\ \text { and } y' \in A(s, C')\}.
\end{multline*} 
\item For all $(\bm{x},s)$ and $n$, $v_n(\bm{x},s,C)$ is convex in $C$. 
\end{enumerate}
\end{proposition}

\proof[\textbf{Proof of Proposition \ref{prop4}}]
\begin{enumerate}[label = (\roman*), wide, labelwidth=!, labelindent=0pt]
\item $y\in A(s,C)$ and $y'\in A(s,C')$ imply $\lambda s \leq \lambda y\leq \lambda(s+C)$ and $(1-\lambda)s \leq (1-\lambda)y' \leq (1-\lambda) (s+C')$; summing terms implies the result.
\item Let $X = (\lambda C + (1-\lambda)C'+ s) $ and $\Delta^S = (X-\xi)/(X-s)$. Note $\Delta^S \in [0,1]$ and $\xi = \Delta^S s + (1-\Delta^S) X$. Let $y=\Delta^S s + (1-\Delta^S)(s+C)$ and $y' = \Delta^S s + (1-\Delta^S)(s+C')$. Then, $y\in A(s, C), y' \in A(s, C')$, and $\lambda y + (1-\lambda)y' = \xi$. 
\item Proof by contradiction follows from  items \ref{item4i} and \ref{item4ii}. 
\item From \autoref{item4iii} and the convexity of $\mathcal{G}_n(\bm{x},y)$ in $y$ for all $n$ and $y$ (by \autoref{prop1} \autoref{item1i}), it follows that 
\begin{multline*}
v_n(\bm{x},s,\lambda C + (1-\lambda)C') 
\\ = \min \{ \mathcal{G}_n(\bm{x}, \lambda y + (1-\lambda)y'): y\in A(s,C), y' \in A(s,C')\}
\\ \leq \min \{ \lambda \mathcal{G}_n(\bm{x},y) + (1-\lambda)\mathcal{G}_n(\bm{x},y'):
\\ y\in A(s,C), y' \in A(s,C')\}
\\ = \lambda v_n(\bm{x},s,C) + (1-\lambda) v_n(\bm{x},s,C').
\end{multline*}
\end{enumerate}
\endproof

Clearly, the assumption that $y_n^*(\bm{x},C) - d_l \leq y_n^*(\bm{\lambda}(\bm{d},\bm{z},\bm{x}),C)$ for all $n$ and all $(\bm{d},\bm{z},\bm{x})$ is in general a challenge to verify \textit{a priori}. Arguments in \cite{Federgruen1986} suggest that as $n$ gets large, $y_n^*(\bm{x},C)$ may converge in some sense to a function $y^*_{\infty}(\bm{x},C)$. From \cite{Malladi2018}, $y_0^*(\bm{x},C)$ is straightforward to determine. Let $\hat{y}(\bm{x},C) \geq y_{\infty}^*(\bm{x},C) \geq y_n^*(\bm{x},C)$ for all $n$ and $\bm{x}$. Then $\hat{y}(\bm{x},C)- d_l \leq y_0^*(\bm{\lambda}(\bm{d},\bm{z},\bm{x}),C)$ for all $(\bm{d},\bm{z},\bm{x})$ implies the above assumption holds. Determination of a function $\hat{y}$ for the general case is a topic for future research. We  present a special case where $y_0^* = y_n^*$ for all $n$ in the appendix section. 

We point out two key differences between the infinite capacity and the finite capacity cases when the reorder cost, $K'=0$. First, when $C$ is infinite, the smallest optimal base stock level $y_n^*(\bm{x})$ is independent of the number of successive approximation steps, making it (relatively) easy to determine. 
Unfortunately, this result may not hold when $C$ is finite except for the situation considered below in  \autoref{prop5}.  This fact has implementation implications for the controllers at the locations; e.g., determining the base stock levels for the capacitated case will in general be more difficult than for the infinite capacity case. 

Second, Propositions \ref{prop3} and \ref{prop4} state that $v(\bm{x}, s, C)$ is non-decreasing and convex in $C$.  We also know that $v(\bm{x}, s, C)$ is convex in $s$ (from \autoref{prop1}, which is also true for the infinite capacity case) and concave and possibly piecewise linear in $\bm{x}$ (from earlier cited results, which is also true for the infinite capacity case).  We showed in \autoref{comput_res} that these structural results can be computationally useful in determining solutions to the stock and production module relocation problem. The relocation problem for determining $(\bm{\Delta^S}, \bm{\sigma}, \bm{u}')$, given $(\bm{x},\bm{s},\bm{u})$, requires knowing $v_l(\bm{x},s_l', u_l')$ for all $ l$. We now consider approaches to compute or approximate $v(\bm{x}, s, C)$, following the presentation of a special case where $y_0^* = y_n^*$ for all $n$. 
\begin{proposition} \label{prop5}
Assume that for all $(d,z,x)$, 
$ y_0^*(\bm{\lambda}(\bm{d},\bm{z},\bm{x}), C) - C \leq y_0^*(\bm{x},C) - d \leq y_0^*(\bm{\lambda}(\bm{d},\bm{z},\bm{x}),C).$
Then, $y_n^*(\bm{x},C) = y_0^*(\bm{x},C)$ for all $n$. 
\end{proposition}

We remark that the left inequality in \autoref{prop5} essentially implies that although capacity may be finite, it is always sufficient to insure the inventory level after replenishment can be $y^*_0(\bm{x},C)$.

\proof[\textbf{Proof of Proposition \ref{prop5}}]
By induction. Assume $y_n^*(\bm{x},C) = y_0^*(\bm{x},C)$. Note therefore, 
$$ v_{n+1}(\bm{x},s,C) = \begin{cases} \mathcal{G}_n(\bm{x},s+C), & s \leq y_0^*(\bm{x},C) - C 
\\ \mathcal{G}_n(\bm{x},s), & s \geq y_0^*(\bm{x},C)
\\ \mathcal{G}_n (\bm{x},y_0^*(\bm{x},C)) & \text{ otherwise.}  \end{cases}$$
Note
\begin{enumerate}[label = (\roman*)]
\item $ \min_y \mathcal{G}_{n+1}(\bm{x},y) \leq \mathcal{G}_{n+1}(\bm{x},y_0^*(x,C))$
\item $ \min_y  \mathcal{G}_{n+1}(\bm{x},y) \geq \min_y \mathcal{L}(\bm{x},y)+ \beta \sum_{d,z} \sigma(d,z,\bm{x}) \min_y v_{n+1} (\bm{\lambda}(d,z,\bm{x}), y-d).$
\end{enumerate}
The minimum with respect to $y$ $v_{n+1}(\bm{\lambda}(d,z,\bm{x}), y-d, C)$ is such that $y_0^*(\bm{\lambda}(d,z,\bm{x}),C) - C \leq y - d \leq  y_0^*(\bm{\lambda}(d,z,\bm{x}),C)$. By assumption, $y = y_0^*(\bm{x},C)$ satisfies these inequalities. Thus, 
\begin{multline*}
\min_y \mathcal{G}_{n+1}(\bm{x},y) \geq \mathcal{L}(\bm{x}, y_0^*(\bm{x},y^*_0(\bm{x},C)) 
\\ + \beta \sum_{d,z} \sigma(d,z,\bm{x}) v_{n+1} (\bm{\lambda}(d,z,\bm{x}), y_0^*(\bm{x},C)-d, C)
\\ = \mathcal{G}_{n+1}(\bm{x},y_0^*(\bm{x},C),
\end{multline*}
and hence $y_{n+1}^*(\bm{x},C) = y_0^*(\bm{x},C)$. 
\endproof
\subsection{Proof of Proposition \ref{prop7} \label{proof_prop7}}
\proof
Let $v_0(\bm{x},s,C ) = \hat{v}_0(\bm{x},s,C) = 0$. 
Consider $\bm{d} =(d_l, \bm{d}_{j\neq l})$, where $\bm{d}_{j \neq l}$ can be considered as additional observation data $z$. Let $\sum_z \sigma(d_l,z,\bm{x}) = \sigma(d_l,\bm{x})$. 
\begin{multline*}
v_1 (\bm{x},s,C)  
\\  = \min_{s\leq y\leq s+C} \bigg\{  \sum_{d_l}\sigma(d_l,\bm{x}) \left [ c(y,d_l) \right ] \bigg\}
\\ = \min_{s\leq y\leq s+C} \bigg\{  \sum_{d_l}\sum_i x_i\sum_j \text{Pr}(j\mid i) \text{Pr}(d_l\mid j)  \left [ c(y,d_l) \right ] \bigg\}
\\ = \min_{s\leq y\leq s+C} \bigg\{  \sum_d\sum_i x_i \bigg( \text{Pr}(d_l\mid i) 
\\ + \sum_j \text{Pr}(j\mid i) \text{Pr}(d_l\mid j) - \text{Pr}(d_l\mid i) \bigg)   \left [ c(y,d_l) \right ] \bigg\}
\\\geq  \min_{s\leq y\leq s+C} \bigg\{  \sum_{d_l}\sum_i x_i \bigg( \text{Pr}(d_l\mid i) -\max_k \text{Pr}(d_l\mid k) 
\\ + \min_k \text{Pr} (d_l\mid k)\bigg)   \left [ c(y,d_l) \right ] \bigg\}
\\\geq  \min_{s\leq y\leq s+C} \bigg\{  \sum_{d_l}\sum_i x_i \text{Pr}(d_l\mid i)  c(y,d_l) \\ - \sum_{d_l}\big(\max_k \text{Pr}(d_l\mid k) - \min_k \text{Pr} (d_l\mid k) \big) c(y,d_l) \bigg\}
\\ \geq  \min_{s\leq y\leq s+C} \bigg\{  \sum_{d_l}\sum_i x_i \text{Pr}(d_l\mid i)  c(y,d_l) \bigg\} 
\\ + \min_{s\leq y\leq s+C} \bigg\{  -\sum_{d_l}\big(\max_k \text{Pr}(d_l\mid k) 
\\ - \min_k \text{Pr} (d_l\mid k) \big) c(y,d_l)  \bigg\}
\\ = \hat{v}_1(\bm{x},s,C) + \min_{s\leq y\leq s+C} \big\{  -\sum_{d_l} k(d_l) c(y,d_l) \big\}
\\ = \hat{v}_1(\bm{x},s,C) - \max_{s\leq y\leq s+C} \big\{ \sum_{d_l} k(d_l) c(y,d_l) \big\}
\\ = \hat{v}_1(\bm{x},s,C) - \sum_{d_l} k(d_l) c(\hat{y},d_l) 
\\ = \hat{v}_1(\bm{x},s,C) - u, \text{ where } u = \sum_{d_l} k(d_l) c(\hat{y},d_l) 
\\ \text{ and } \hat{y} \in \{s,s+C\} 
\text{ due to convexity of } c({y},d_l)
\\ \forall \ y, d_l, \text{ where } k(d_l) = \big(\max_k \text{Pr}(d_l\mid k) - \min_k \text{Pr} (d_l\mid k) \big). 
\end{multline*}
By induction and infinite summation, 
\begin{multline*} v_n(\bm{x},s,C) \geq \hat{v}_n(\bm{x},s,C) - u(1+\beta +\dots +\beta^n);
\\ v(\bm{x},s,C) \geq \hat{v}(\bm{x},s,C) - u/(1-\beta) .\end{multline*}
\endproof
\subsection{The Heuristic LARRO \label{ex_hr}}
We now present a heuristic for large instances with low computational overhead.  LARRO stands for lookahead of rollout for relocations only. 
\begin{align}
 &\nonumber \text{LARRO:} 
 \\ & \nonumber \min_{\bm{\Delta^S}, \bm{u}',\bm{y}} \sum_l  \bigg\{ (K^{S+}_l \Delta^{S+}_l + K^{S-}_l \Delta^{S-}_l ) 
 \\ \nonumber & + K^M\sum_l \abs{u_l-u_l'}/2 + (\zeta_l + \eta_l)/2  \bigg\},
\\ \nonumber & \text{subject to }  
\\ & \zeta_l \geq \gamma_j^l (s_l + \Delta^{S+}_l - \Delta^{S-}_l) + \hat{\gamma}_j^l \ \forall \ (\gamma_j^l, \hat{\gamma}_j^l) \in \Gamma^l_{t+1}(u_l)  \ \forall \ l   \nonumber
\\ & \eta_l  \geq \theta_j^l u_l' + \hat{\theta}_j^l  \ \ \ \ \forall \ ( \theta_j^l, \hat{\theta}_j^l) \in \Theta^l_{t+1}(s_l)  \ \forall \ l \nonumber
\\ \nonumber & \sum_l u_l' = Y
 \\ \nonumber & \sum_l \Delta^{S+}_l = \sum_l \Delta^{S-}_l,
  \\ \nonumber & 0\leq u_l' \leq Y_l', \ \forall \ l
  \\ \nonumber & 0\leq \Delta^{S+}_l \leq \sum_{k\neq l} (s_k)^+, \ \forall \ l
  \\ \nonumber & 0\leq \Delta^{S-}_l \leq - (s_l)^+, \ \forall \ l
 \\ & u_l',\ \Delta^{S+}_l, \Delta^{S-}_l \in \mathbb{Z}, \ \ \eta_l, \ \zeta_l \in \mathbb{R}  \ \ \forall \ l
 \end{align} 

\begin{proposition} LARRO can be solved exactly by relaxing the integrality constraints. 
\end{proposition}
\subsection{\label{ex_tables} Results: Additional Tables}
We now present additional numerical results that complement \autoref{comput_res}.
\begin{table}[h!]
\centering
\caption{Variation of average savings due to \ac{RRO} over \ac{DNF} with varying $\theta$  for Instance Set A}
\label{glr-blend-a}
\begin{tabular}{cccccc}
\toprule
$G \ \backslash \ \theta$       & 0    & 0.2  & 0.5  & 0.8  & 1    \\
\hline
1       & -131\% & 42\% & 40\% & 39\% & 38\% \\
2       & -133\% & 38\% & 36\% & 35\% & 35\% \\
5       & -141\% & 44\% & 42\% & 40\% & 39\% \\
\hline
Overall & -135\% & 41\% & 39\% & 38\% & 37\%
\\\bottomrule  \hline
\end{tabular}
\end{table}

\begin{table}[h!]
\centering
\caption{Variation of average savings due to \ac{LSF} over \ac{DNF} across $\theta$  for Instance Set A}
\label{laj-blend-a}
\begin{tabular}{cccccc}
\toprule
$G \ \backslash \ \theta$       & 0    & 0.2  & 0.5  & 0.8  & 1    \\
\hline
1       & 34\% & 43\% & 37\% & 33\% & 31\% \\
2       & 38\% & 38\% & 33\% & 29\% & 27\% \\
5       & 44\% & 44\% & 37\% & 33\% & 31\% \\
\hline
Overall & 39\% & 42\% & 36\% & 32\% & 30\% 
\\\bottomrule  \hline
\end{tabular}
\end{table}

\begin{table}[h!]
\centering
\caption{Variation of average savings due to \ac{RSF} over \ac{DNF} across $\theta$  for Instance Set A}
\label{jr-blend-a}
\begin{tabular}{cccccc}
\toprule
$G \ \backslash \ \theta$       & 0    & 0.2  & 0.5  & 0.8  & 1    \\
\hline
1       & 37\% & 44\% & 39\% & 35\% & 32\%
\\\bottomrule  \hline
\end{tabular}
\end{table}

\begin{table}[h!]
\centering
\caption{Variation of average savings due to heuristics over \ac{DNF} across $G$ for $\theta = 0.2$ on a shorter horizon $T=10$ instead of $T=30$  for Instance Set A}
\label{all_seta_0.2_shorter_hor}
\begin{tabular}{ccccc|cc}
\toprule
$G$     & \ac{MNF}  & \ac{MP} & \ac{RRO} & \ac{LSF} & \ac{LSF}-SS & \ac{LSF}-CO \\
\hline
1       & -3\% & -2\% & 24\% & 28\% & 26\% & 28\% \\
2       & -4\% & 6\%  & 22\% & 24\% & 23\% & 25\% \\
5       & -5\% & 14\% & 27\% & 29\% & 27\% & 29\% \\
\hline
Overall & -4\% & 6\%  & 24\% & 27\% & 25\% & 28\%
\\\bottomrule  \hline
\end{tabular}
\end{table}

\begin{table}[h!]
\centering
\caption{Value of mobility (\% savings over \ac{DNF}) using \ac{LSF} with $\theta = 0.2$ across varying $ K^S$ and $K^M$ for $G=1$ instances of Instance Set A}
\label{val_laj_seta_g1}
\begin{tabular}{cc|ccccc}
\toprule
                         &     & \multicolumn{5}{c}{Module movement cost $K^M$}  \\
                           &  &0   & 1.5      & 2        & 2.5 &1000       \\ 
\hline
&0    & 53\% & 50\% & 49\% & 52\% & 55\% \\
&1.5  & 50\% & 44\% & 38\% & 42\% & 46\% \\
&2    & 54\% & 42\% & 44\% & 42\% & 43\% \\
&2.5  & 48\% & 44\% & 40\% & 36\% & 40\% \\
\multirow{-6}{*}{\rotatebox[origin=c]{90}{\parbox{2.3cm}{\centering  \scriptsize{Transshipment} \\  cost $K^S$}}} &1000 & 52\% & 43\% & 45\% & 38\% & -3\%  
\\ \bottomrule  \hline
\end{tabular}
\end{table}

\begin{table}[h!]
\centering
\caption{Value of mobility (\% savings over \ac{DNF}) using \ac{LSF} with $\theta = 0.2$ across varying $ K^S$ and $K^M$ for $G=5$ instances of Instance Set A}
\label{val_laj_seta_g5}
\begin{tabular}{cc|ccccc}
\toprule
                         &     & \multicolumn{5}{c}{Module movement cost $K^M$}  \\
                           &  &0   & 1.5      & 2        & 2.5 &1000       \\ 
\hline
&0    & 48\% & 49\% & 47\% & 53\% & 50\% \\
&1.5  & 49\% & 45\% & 44\% & 46\% & 39\% \\
&2    & 50\% & 44\% & 41\% & 44\% & 37\% \\
&2.5  & 49\% & 46\% & 45\% & 45\% & 34\% \\
\multirow{-6}{*}{\rotatebox[origin=c]{90}{\parbox{2.3cm}{\centering  \scriptsize{Transshipment} \\  cost $K^S$}}}  &1000 & 50\% & 48\% & 50\% & 44\% & -5\% 
\\ \bottomrule  \hline
\end{tabular}
\end{table}
\end{document}